\documentclass[12pt]{article}
\usepackage[cp1251]{inputenc}
\usepackage{mathrsfs}                                    
\usepackage[T2A]{fontenc}           

\usepackage[english,russian]{babel} 

\usepackage{amssymb}
\usepackage{amsmath}

\newtheorem{theorem}{Theorem}[section]

\newtheorem{lemma}{Lemma}[section]

\newtheorem{corollary}{Corollary}[section]
\numberwithin{equation}{section}

\input cyracc.def

\begin{document}
\title{The Sharp Lower Bound  \\of Asymptotic Efficiency of Estimators \\ in the Zone  of Moderate Deviation Probabilities} \author{
Mikhail Ermakov\\Mechanical Engineering Problem Institute RAS.\\
Bolshoy pr. V.O. 61\\ St.Petersburg 199178 \\RUSSIA. }
\maketitle
\centerline{\it e-mail: erm2512@mail.ru}

\begin{abstract}  For the zone of moderate deviation probabilities the local asymptotic minimax lower bound of asymptotic efficiency
of estimators is established.
The estimation parameter is multidimensional. The lower bound admits the
interpretation as the lower bound of asymptotic efficiency in confidence estimation.
\end{abstract}

\newpage

 \section{Introduction}

The asymptotic normality of estimators is  a key property allowing to construct confidence sets if the sample size is sufficiently large.
The problem of accuracy of the normal approximation   emerges simultaneously with its implementation.
   The inequalities of the Berry-Esseen type  and the Edgeworth expansions (see \cite{fe,bar,re,gu} and references therein)
show that the convergence rate to the  normal distribution
has the order
$n^{-1/2}$ ( here $n$ is a sample size).  The significant levels $\alpha$ of confidence sets have usually  small values
( $\alpha = 0.1; 0.05; 0.01$ are the standard values in practice
).   For such  small values of $\alpha$ the rate of convergence $n^{-1/2}$ does not allow to talk
about adequate accuracy of approximation  for the sample sizes of several hundreds observations  or smaller.
 From this viewpoint the study of  asymptotic properties of estimators
in the zones of  large and moderate deviation probabilities is of special interest.
The problem of  lower bounds
for asymptotic efficiency in
these zones emerges as well.  The  asymptotic efficiency of estimators
in the zone of large deviation probabilities is  analyzed on the base of Bahadur efficiency \cite{ba, van, ra, pu}.

The study of large deviation probabilities of estimators is a
rather difficult problem. This problem is often replaced with the study of their moderate deviation probabilities.
Let $X_1,\ldots,X_n$ be independent sample  of random variable $X$
having the probability measure $P_\theta, \theta \in R^1$. Let
$b_n > 0,  b_n \to 0, nb_n^2 \to \infty$ as $n \to \infty$. Let
 $\theta_0 \in R^1$. Then (see \cite{er03}) for any estimator
$\hat \theta_n$
\begin{equation}\label{1.1}
\lim\inf_{n\to\infty} \inf_{\theta = \theta_0,\theta_0 + 2b_n}
(\frac {1}{2}nb_n^2)^{-1} \ln P_\theta( |\hat \theta_n - \theta| >
b_n) \ge - I(\theta_0).
\end{equation}
Here we suppose that there exists the finite Fisher information  $I(\theta)$ for all $\theta$ in some vicinity of
$\theta_0$. Note that the lower bound of the local Bahadur asymptotic efficiency is a particular case of  (\ref{1.1}).

The natural problem arises on the quality of logarithmic approximation  for the obtaining
 confidence sets.
The  distributions of estimators admit usually the approximation by the  sums
 $\bar X =
n^{-1}(X_1+\ldots+X_n)$ of independent random variables.
(see \cite{se,van,gu} and references therein).    Thus it is of interest to compare for the sample mean
 $\bar X$ the confidence intervals obtained by the normal approximation and the basic term of logarithmic asymptotic.
If there exists an exponential moment
 $E [\exp\{t|X_1|\}] < C <\infty, t>0$, the sample mean $\bar X$ satisfies the Bernstein inequality
\begin{equation}\label{1.2}
P(n^{1/2}(\bar X - E [X_1]) > x) < \exp \left\{-\frac{x^2}{2\sigma^2}(1 +
o(1))\right\}, \quad x> x_0
\end{equation}
with $\sigma^2 = \mbox{Var} [X_1]$.

The confidence interval based on the main term of asymptotics of right-hand side of (\ref{1.2}) is the following
\begin{equation}\label{1.3}
\left(\bar X - \frac{\sigma\sqrt{2|\ln(\alpha/2)|}}{\sqrt{n}},
\bar X+ \frac{\sigma\sqrt{2|\ln(\alpha/2)|}}{\sqrt{n}}\right)
\end{equation}
 instead of the standard one
\begin{equation}\label{1.4}
\left(\bar X - x_{\alpha/2}\frac{\sigma}{\sqrt{n}}, \bar X +
x_{\alpha/2}\frac{\sigma}{\sqrt{n}}\right)
\end{equation}
where $x_{\alpha/2}$ satisfies $\alpha/2= \Phi(-x_{\alpha/2})$. Here $\Phi(x)$ is the standard normal distribution
function.

If $\alpha = 0.1; 0.05; 0.01$ respectively the confidence intervals defined by
(\ref{1.3}) are the following
$$
(\bar X - 2.44\frac{\sigma}{\sqrt{n}}, \bar X +
2.44\frac{\sigma}{\sqrt{n}}),
$$
$$
(\bar X - 2.71\frac{\sigma}{\sqrt{n}}, \bar X +
2.71\frac{\sigma}{\sqrt{n}}),
$$
$$
(\bar X - 3.25\frac{\sigma}{\sqrt{n}}, \bar X +
3.25\frac{\sigma}{\sqrt{n}})
$$
instead of the standard ones defined by the normal approximation (\ref{1.4})
$$
(\bar X - 1.65\frac{\sigma}{\sqrt{n}}, \bar X +
1.65\frac{\sigma}{\sqrt{n}}),
$$
$$
(\bar X - 1.96\frac{\sigma}{\sqrt{n}}, \bar X +
1.96\frac{\sigma}{\sqrt{n}}),
$$
$$
(\bar X - 2.576\frac{\sigma}{\sqrt{n}}, \bar X +
2.576\frac{\sigma}{\sqrt{n}}).
$$
If $\alpha = 0.1; 0.05$, the implementation of (\ref{1.3})
requires the doubling of the  number of observations for obtaining the same
width of confidence interval as in  (\ref{1.4}). At the same time the normal approximation works
in a rather narrow zone of moderate deviation probabilities in comparison with  the Bernstein inequality (\ref{1.2}).
 Thus the analysis of confidence intervals on the base of logarithmic asymptotics of large and moderate deviation
 probabilities is also reasonable. It should be noted that
there exist powerful methods for constructing
 accurate boundaries of confidence intervals such as asymptotic expansions
    (see \cite{fe, gu, bar,  re, sk} and references therein), bootstrap  (see  \cite{ef,bu,van,gu}
and references therein) and so on.

For the zone of moderate deviation probabilities the normal approximation of statistics is the subject of
numerous publications (see \cite {bar, al, bu, gu, re, in, ju} and references therein).
The goal of the paper is to prove the sharp local asymptotic minimax  lower bound for the estimators in this zone.
The estimation parameter is multidimensional.
For one - dimensional parameter the local asymptotic minimax  lower bound for the sharp asymptotics
of moderate deviation probabilities  of estimators
has been  established in \cite{er03}. Thus the local asymptotic minimax lower bound for estimators \cite{ha,ib,le, str, van}
is extended on the zone of moderate deviation
probabilities.

We make use of the letters C and c as generic notation for positive constants.
Denote $\chi(A)$ the indicator of set $A$, $[a]$ - the integral part of $a$.
For any $u,v \in R^d$   denote  $u'v$ the inner product  of $u, v$ and
$u'$ the transposed vector of $u$. For positive sequences $a_n$ denote $a_n \asymp b_n$, if $c < a_n/b_n < C$, and denote
$a_n >>> b_n$ if $a_n/b_n \to \infty$ as $n \to \infty$. For any set of events $B_{...}$ denote $A_{...}$
the complementary event to $B_{...}$.

\section{ Main Result}

Let $X_1,\ldots,X_n$ be i.i.d.r.v.'s having a probability measure (p.m.)  $P_\theta, \theta \in \Theta \subseteq R^d$, defined on a probability space
 $(S,\Upsilon)$.
Suppose  p.m.'s $P_\theta, \theta \in \Theta$,
are absolutely continuous w.r.t. p.m.  $\nu$ defined on the same probability space $(S,\Upsilon)$. Denote $f(x,\theta)
= \frac{dP_\theta}{d\nu}(x), x \in S$.  For any $\theta_1, \theta_2 \in R^d$  denote
$P^a_{\theta_1,\theta_2}$ and $P_{\theta_1,\theta_2}^s$ respectively absolutely continuous and singular components
of p.m. $P_{\theta_1}$ w.r.t. $P_{\theta_2}$. For all $x
\in S$ such that $f(x,\theta) \ne 0$  denote  $ g(x,\theta,\theta+u) =
(f(x,\theta+u)/f(x,\theta))^{1/2} - 1, u \in R^d$.

The statistical experiment  $\Psi= \{(S,\Upsilon), P_\theta, \theta \in R^d\}$ has the finite Fisher information at the point
 $\theta \in R^d$ if there exists the vector function
  $\phi_\theta(x) = (\phi_{\theta,1}(x),\ldots,
 \phi_{\theta,d}(x))', x \in S , \phi_{\theta,i}\in L_2(P_\theta), 1 \le i \le d$ such that
$$
\int_S \left(g(x,\theta,\theta+u) - \frac{1}{2}u'\phi_\theta(x)\right)^2 dP_\theta = o(|u|^2), \quad
P^s_{\theta+u,\theta}(S)  = o(|u|^2)
$$
as $u \to 0$.

The Fisher information matrix at the point
  $\theta$ equals
$$
I(\theta) = \int_S \phi_\theta \phi'_\theta\,\, dP_\theta .
$$
For any $P_{\theta_1}, P_{\theta_2}, \theta_1,\theta_2 \in R^d$
the Hellinger distance equals
$$
\rho(P_{\theta_1},P_{\theta_2}) = \rho(\theta_1,\theta_2) =
\left(\int_S (f^{1/2}(x,\theta_1) - f^{1/2}(x,\theta_2))^2 \,d\nu\right)^{1/2}.
$$
We make the following assumptions.

Let $\theta_0 \in \Theta$ and let  $\Theta$ be open set. Let $0 < \lambda \le 1$.

\noindent{\bf A1.} For all $\theta$  in some vicinity $\Theta_0$ of   the point  $\theta_0 \in \Theta$ there exists  the positive definite Fisher information matrix
$I(\theta)$.

\noindent{\bf A2.} For all
 $\theta, \theta + u \in \Theta_0$   there hold
\begin{equation}\label{2.1}
\int_S (g(x,\theta,\theta+u) - \frac {1}{2}u'\phi_\theta(x))^2 \,dP_\theta < C|u|^{2+\lambda},\quad
P^s_{\theta+u,\theta}(S) < C|u|^{2+\lambda},
\end{equation}
\begin{equation}\label{2.2}
|4\rho^2(\theta,\theta+u) - u'I(\theta)u| < C|u|^{2+\lambda},
\end{equation}
\begin{equation}\label{2.3}
\int_S |\phi_\theta(x)|^{2 + \lambda}\, dP_\theta < C<\infty,
\end{equation}
\begin{equation}\label{2.4}
h'I(\theta)h - h'I(\theta+u)h < C |h|^2|u|^\lambda.
\end{equation}
The constants $C$ in  (\ref{2.1}-\ref{2.4}) do not depend on  $\theta, \theta + u \in \Theta_0$.

We say that a set $\Omega \subset R^d$ is central-symmetric if
 $x \in \Omega$ implies $-x \in \Omega$.

We make the following assumptions

\noindent{\bf B1.} The set $\Omega$ is convex and central-symmetric.

\noindent{\bf B2.} The boundary  $\partial\Omega$ of the set $\Omega$ is $C^2$-manifold.

\noindent{\bf B3.} The principal curvatures at each point of $\partial\Omega$ are negative.

Denote $\zeta$- Gaussian random vector in $R^d$ such that
$E\zeta = 0, E[\zeta\zeta'] = I$. Here $I$ is the unit matrix.

\begin{theorem} \label{t1} Assume A1, A2 and B1-B3. Let $nb_n^2 \to
\infty, nb_n^{2+\lambda} \to 0$, $b_n-b_{n-1} = o(n^{-1}b_n^{-1})$ as $n \to \infty$. Then for any estimator
$\hat\theta_n = \hat\theta_n(X_1,\ldots,X_n)$
\begin{equation}\label{2.7}
\liminf_{n\to\infty}\sup_{|\theta-\theta_0|<C_nb_n}\frac{P_\theta
(I^{1/2}(\theta_0)(\hat\theta_n
- \theta) \notin b_n\Omega)} {P(\zeta\notin
n^{1/2}b_n\Omega)}\ge 1
\end{equation}
with $C_n \to \infty$ as $n \to \infty$.
\end{theorem}

Wolfowitz \cite{wo} was the first who pointed out the relationship of lower bounds of
(\ref{2.7})-type
with the problem of asymptotic efficiency in the confidence estimation.

In \cite{er03} Theorem \ref{t1}  has been established for $\theta \in \Theta \subseteq R^1$ if (\ref{2.1})-(\ref{2.3}) is valid.
If $d = 1$,  (\ref{2.4}) follows from (\ref{2.2}).
Note that (\ref{2.4}) is fulfilled evidently in the case of location  parameter.
  If (\ref{2.4}) does not valid, we could not  take $I^{1/2}(\theta_0)$ as the constant normalized matrix  in
 (\ref{2.7}).

In confidence estimation the set $\Omega$ is usually a   ball $\Omega_r$ having the center zero and the radius $r>0$.
In this case the asymptotic of denominator in (\ref{2.7}) is known.

\begin{corollary} Let assumptions of Theorem \ref{t1} be valid. Let $\Omega = \Omega_r$.Then for any estimator
$\hat\theta_n = \hat\theta_n(X_1,\ldots,X_n)$
\begin{equation}\label{2.8}
\liminf_{n\to\infty}\sup_{|\theta-\theta_0|<C_nb_n} 2^{d/2 - 1}\Gamma(d/2) (n^{1/2} b_n r)^{2-d} \exp\{nb_n^2r^2/2\}   P_\theta
(I^{1/2}(\theta_0)(\hat\theta_n
- \theta) \notin b_n\Omega_r) \ge 1
\end{equation}
with $C_n \to \infty$ as $n \to \infty$.
\end{corollary}

If $\Omega$ is the ellipsoid $\Omega_{\sigma,r} = \left\{\theta:  \sum_{i=1}^d \sigma_i^2\theta_i^2 >r^2,
\theta = \{\theta_i\}_{i=1}^d, \theta_i \in R^1\right\},
\sigma = \{\sigma_i\}_{i=1}^d, \sigma_1=\sigma_2=\ldots=\sigma_k > \sigma_{k+1} > \ldots > \sigma_d  > 0,$ we get the following asymptotic
(see \cite{li}) in the denominator of (\ref{2.7})
\begin{equation}\label{2.9}
P(\zeta\notin
n^{1/2}b_n\Omega_{\sigma,r}) =  C_k (n^{1/2}b_nr)^{k-2}
\exp\{-nb_n^2r^2/2\}(1+o(1)).
\end{equation}
Here $C_k = 2^{1-k/2}\sigma_1^{1-k}(\Gamma(k/2))^{-1}\prod_{i=k+1}^d (1 - \sigma_r^2/\sigma_1^2)^{-1/2}$.

The assumptions of Theorem \ref{t1} are rather weak. The sharp asymptotics of moderate deviation
probabilities of likelihood ratio were established
 under the more restrictive assumptions (see
\cite{bar,bor,bu,sk} and references therein).  The lower bounds for moderate deviation probabilities do not
require such strong assumptions (see \cite{ar, er03}) and are usually
 proved more easily than the
upper bounds.

The assumptions of Theorem \ref{t1}
  are different from the traditional assumption of local asymptotic normality.
Thus Theorem \ref{t1} could not be straightforwardly extended on the models having this property.
 At the same time  A1,A2 represent slightly more stable form of
usual assumptions arising in the proof of local asymptotic normality. This allows to make use of the technique
arising in the proofs of local asymptotic normality and to get the results similar to (\ref{2.7}) for other models of estimation.
This problem will be considered in the sequel.

For the semiparametric estimation the local asymptotic minimax lower bounds
in the zone of moderate deviation probabilities have been established in \cite{er04}. In \cite{er04}
the statistical functionals take the values in $R^1$. The results were based on the assumptions that (\ref{2.1}-\ref{2.3}) hold uniformly
for the families of "least-favourable" distributions.
In the case of multidimensional parameter  the additional assumptions (\ref{2.4}) arises only.
 Thus the difference is not very significant.

The plan of the proof of Theorem \ref{t1} is the following. In section 3 we outline the basic steps of the proof. After that the proof
 are given for the most simple geometry of the set $\Omega$.
For the  arbitrary geometry of set $\Omega$ we point out the differences in the proof
at the end of section 3. The key  Lemmas \ref{l3.3}, \ref{l3.4} are proved in section 4.
 The proof of Lemma \ref{l3.4} is based on new
 Theorems \ref{t2} and \ref{t3} on large deviation probabilities of sums of independent random vectors.
The proofs of Theorems \ref{t2} and \ref{t3}
are given in section 5. The proofs of  technical Lemmas of sections 3 and 4 are given in section 6.
\section{Proof of Theorem \ref{t1}}

To simplify the notation we suppose that  $\theta_0$ equals zero.
Suppose the matrix $I(\theta_0)$ is the unit.

For any $\theta_1,\theta_2 \in \Theta$ denote
$$
\xi_s(\theta_1,\theta_2) =
\ln\frac{f(X_s,\theta_2)}{f(X_s,\theta_1)},\quad
\tau_s(\theta_1) = \{\tau_{ks}(\theta_1)\}_1^d =\phi_{\theta_1}(X_s)
$$
with $1\le s \le n$.

We will often omit
  $\theta=\theta_0$ in notation. For example, we shall write $\xi_s(\theta) = \xi_s(\theta_0,\theta),   \tau_s =
\tau_s(\theta_0)$. The index
 $s$ will be omitted for $s=1$.
For example, $\tau = \tau_1(\theta_0)$.

Denote $\psi_n=
n^{-1/2}I^{-1/2}(\theta_0)\sum_{s=1}^n \tau_s$. Note, that
$(\theta-\theta_0)'\sum_{s=1}^n \tau_s$ is the stochastic part of the
linear approximation  of logarithm of likelihood ratio.

The reasoning is based on the standard proof of local asymptotic minimax lower bound
 \cite{ha,ib,le, str, van}. In particular we make use of the fact that the minimax risk
exceeds the Bayes one and study the asymptotic  of Bayes risks.
However, in this setup, the
estimates of  residual terms of  asymptotics of posterior Bayes risks should have the order
$o(\exp\{-cnb_n^2\})$. This does not allow to make use of the technique of local asymptotic normality
\begin{equation}\label{3.1}
\sum_{s=1}^n \xi_s(u_n) - n^{1/2}u'_nI^{1/2}\psi_n + \frac{1}{2}nu_n' I u_n = o_P(1)
\end{equation}
in the zone $|u_n| \le Cb_n$ of moderate deviation probabilities. This is the basic reason of differences in the proof.

Instead of  (\ref{3.1}) we are compelled to prove that, for any $\epsilon > 0$,
\begin{equation}\label{3.2}
P\left(\sup_{u \in U_n}\left\{\sum_{s=1}^n \xi_s(u) -
n^{1/2}u' I^{1/2}\psi_n + \frac{1}{2}nu' I u\right\} >
\epsilon\right) = o(\exp\{-cnb_n^2\})
\end{equation}
where $U_n$ is a fairly broad set of
 parameters.
Therefore, the main problem is how to narrow down the set
 $U_n$.

The following two facts have allowed to solve this problem.

The normalized values of
posterior Bayes risks tend to a constant
in probability.

In the zone of moderate deviation probabilities the normal approximation   \cite{bahr,os}  holds
for the sets of events
$ \psi_n \in n^{1/2}\Gamma_{ni}$  where the domain
 $\Gamma_{ni}$ has a diameter $o(n^{-1}b_n^{-1})$.

Thus we can find the asymptotic of posterior Bayes risk independently for
 each an event
$\psi_n \in n^{1/2}\Gamma_{ni} $ , sum over $i$  and get the lower bound.
Fixing the set $\Gamma_{ni}$
 allows to replace the proof of  (\ref{3.2}) with
\begin{equation}\label{3.2a}
\begin{split}&
P\left(\sup_{u \in U_n}\left\{\sum_{s=1}^n \xi_s(u) -
n^{1/2}u'I^{1/2}\psi_n + \frac{1}{2}nu' I u, \right\} >
\epsilon, \psi_n \in n^{1/2}\Gamma_{ni}, A_{1n}\right)\\& =
o\left(\int_{n^{1/2}\Gamma_{ni}}\exp\{-x^2/2\} dx\right)
\end{split}
\end{equation}
where $P(A_{1n}) = 1+o(1)$.

To narrow down the sets $U_n $ we  define
  the lattice $\Lambda_n$ in the cube $K_{v_n}, v_n = Cb_n $ and
 split $\Lambda_n$ into subsets $\Lambda_{nile}$. The set $ \Lambda_ {nile} $
 is the lattice in the union of a finite number of very narrow
 parallelepipeds $K_{nij}$ whose orientation  is given by the
 position of the set $ \Gamma_ {ni} $ relative to $ \theta_0 $.  The problem of
Bayes risk minimization is solved independently for each set $\Lambda_ {nile} $ and the results are added.

 Note that the proof of
 (\ref{3.2a}) with $U_n= \Lambda_{nile}$ is based on the  "chaining method" together with the inequality
\begin{equation}\label{3.2b}
\begin{split}&
P\left(\sum_{s=1}^n \xi_s(\theta_1,\theta_2) -
(\theta_2-\theta_1)' \sum_{s=1}^n \tau_{s\theta_1}
+ \frac{1}{2}n(\theta_2-\theta_1)' I (\theta_2 -\theta_1) >
\epsilon, \right. \\& \left. \psi_n \in n^{1/2}\Gamma_{ni}, A_{1n}\right) \le C
|\theta_2-\theta_1|^2 b_n^\lambda\int_{n^{1/2}\Gamma_{ni}}\exp\{-x^2/2\} dx.
\end{split}
\end{equation}
To prove (\ref{3.2b}) we implement simultaneously the Chebyshev inequality to the first sum in the left-hand side of (\ref{3.2b})
and theorem on large deviation probabilities for
$\psi_n$. Thus we
 prove some anisotropic version of  theorem on large deviation probabilities
  (see Theorem \ref{t3}).

 Denote $v_n = Cb_n$. Define a sequence
  $\delta_{1n}=
c_{1n}(nb_n)^{-1}$, with $c_{1n}\to 0, c_{1n}^{-3}nb_n^{2+\lambda}\to 0$
as $n \to \infty$. In the cube $K_{v_n}=[-v_n,v_n]^d$ we define a lattice
 $ \Lambda_n = \{h: h= (j_1\delta_{1n},\ldots,j_d\delta_{1n}), -l_n \le j_k \le l_n=
[v_n/\delta_{1n}], 1 \le k \le d\}$. Thus
$l_n\asymp c^{-1}_{1n}nb_n^2$.

We split the cube  $K_{\kappa v_n}, 0< \kappa <1$ on the small cubes $\Gamma_{ni} = x_{ni} +
(-c_{2n}\delta_{1n}, c_{2n}\delta_{1n}]^d$, where $c_{2n} \to
\infty, c_{2n}\delta_{1n}= o(n^{-1} b_n^{-1}), c_{2n}^3c_{1n}^{-3}nb_n^{2+\lambda}\to 0$ as $n \to \infty,
1 \le i \le m_n=[(\kappa c_{2n}^{-1}C c_{1n}^{-1})^dn^db_n^{2d}], x_{ni} \in K_{v_n}$.

Suppose $C$ is chosen so that $b_n\Omega \subset K_{(1-\kappa)v_n}$.

For each $x_{ni}, 1 \le i
\le m_n$ we define the partition of the cube $K_{v_n}$ on the subsets
$$K_{nij}
=K(\theta_{nij})=\{x: x = \lambda x_{ni} + u +\theta_{nij}, u =
\{u_k\}_{k=1}^d,
$$
$$
 u \bot x_{ni}, |u_k| \le c_{3n} \delta_{1n},
\lambda \in R^1, u \in R^d\} \cap K_{v_n}, 1 \le j \le m_{1ni},
$$
where $c_{3n}/c_{2n} \to \infty, c_{3n}\delta_{1n}= o(n^{-1}
b_n^{-1}), c_{3n}^3c_{1n}^{-3}nb_n^{2+\lambda}\to 0$ as $n \to \infty$.

Let us fix   $ i $.
 Suppose  $x_{ni}$ is parallel to  $e_1=(1,0,\ldots,0)'$.
This does not cause serious differences in the reasoning. Denote $\Pi_1$
the subspace orthogonal to $e_1$. Suppose the points
 $ \theta_ {nij}, 1 \le j \le m_{1ni} $ are chosen so that they form a lattice in $ \Pi_1 \cap
K_{v_n}$.  Define the sets
$$
\Lambda_n(\theta_{nij}) =
K(\theta_{nij}) \cap \Lambda_n, 1 \le j \le m_{1ni}, \quad \Theta_{ni} =
\{\theta: \theta= \theta_{nij}, 1 \le j \le m_{1ni}\}.
$$

The risk asymptotic is defined by the set
\begin{equation}\label{3.4}
M = \{ x : |x| = \inf_{y\in \partial \Omega} |y|,\quad x \in
\partial \Omega\,\}.
\end{equation}
We begin with the proof of Theorem \ref{t1}  for the two-point case $M=\{-y,y\}, y \in \partial \Omega$.
For arbitrary geometry of the set $M$ we are compelled to make use of a rather cumbersome constructions.
At the same time the basic part of the proof is the same.

Let $\theta_{nij_0}$ be such that $b_ny \in K(\theta_{nij_0})$  Then $-b_ny \in K(-\theta_{nij_0})$. Let us split $\Theta_{ni}$ on the subsets
\begin{equation}\label{3.a6a}
\begin{split}&
\Theta_{i}(k_1,\ldots,k_{d-d_1})= \{\theta: \theta=
\theta_{nij_0}  + (-1)^{t_2} 2k_2 c_{3n}\delta_{1n}e_{2}\\& +
\ldots + (-1)^{t_{d}} 2k_{d}c_{3n}\delta_{1n} e_d;\,\,
t_2,\ldots t_{d}= \pm 1\}
\end{split}
\end{equation}
where $0 \le k_2,\ldots,k_{d} < C_{1n}$ with  $C_{1n}c_{3n}c_{1n}  \to \infty, nC_{1n}^{3}c_{3n}^{3}c_{1n}^{3}b_n^{2+\lambda} \to 0$ as $n \to \infty$.

Denote
\begin{equation}
\tilde K_{ni}(k_1,\ldots,k_{d-1}) = \cup _{\theta \in \Theta_i(k_1,\ldots,k_{d-1})} K(\theta).
\end{equation}
It will be convenient to number the sets $ \tilde
K _ {ni} (k_1, \ldots, k _ {d-d_1}) $ denoting their $ \tilde K _ {ni1},
\ldots,  \tilde K_{nim_{2ni}}$. Denote
\begin{equation}\label{3.7a}
\Theta_{nie}=\Theta_{ni}\cap \tilde K_{nie}, \quad \Lambda_{nie} = \tilde
K_{nie}\cap \Lambda_n, \quad 1 \le e \le m_{2ni}.
\end{equation}
Thus $\Theta_{nie}$ contains $k= 2^{d-1}$ points, that is, $\Theta_{nie}= \{\theta_{j}\}_{j=1}^{k}$.

In this notation the problem of risk minimization on $\Lambda_n$ is reduced to the same problems  on the subsets
$\Lambda_{nie}$
\begin{equation}
\label{3.8}
\begin{split}&
\inf_{\hat\theta_n}\sup_{\theta \in K_{v_n}} P_\theta(\hat\theta_n
- \theta \notin b_n\Omega)\\& \ge \inf_{\hat\theta_n}\,
(2l_n)^{-d}\sum_{i=1}^{m_n}\sum_{\theta \in \Lambda_n}  P_\theta(\hat\theta_n -\theta \notin b_n \Omega,
 \psi_n \in n^{1/2}\Gamma_{ni})\\&\ge
(2l_n)^{-d}\sum_{i=1}^{m_n}\sum_{e=1}^{m_{2ni}}\inf_{\hat\theta_n}
\sum_{\theta \in \Lambda_{nie}}  P_\theta(\hat\theta_n - \theta
\notin b_n \Omega, \psi_n \in n^{1/2}\Gamma_{ni}).
\end{split}
\end{equation}
Thus we can minimize the Bayes risk  on each  subset  $\Lambda_{nie}$ independently and  make use
of the own linear approximation (\ref{3.1}) of logarithms of likelihood ratio
on each set $U_n =\Lambda_{nie}$.

For the arbitrary geometry of the set $M$ the additional summation over index $l, 1 \le l \le m_{3ni}$
caused the different points of $M$ arises in (\ref{3.8}). Thus the right-hand side of (\ref{3.8}) is the following
\begin{equation}
\label{3.8a}
(2l_n)^{-d}\sum_{i=1}^{m_n}\sum_{l=1}^{m_{3ni}}\sum_{e=1}^{m_{2nil}}\inf_{\hat\theta_n}
\sum_{\theta \in \Lambda_{nile}}  P_\theta(\hat\theta_n - \theta
\notin b_n \Omega, \psi_n \in n^{1/2}\Gamma_{ni}).
\end{equation}
The definition of the sets $\Lambda_{nile}$ is akin to $\Lambda_{nie}$. The statement (\ref{3.8}) with the right-hand
side (\ref{3.8a}) is the  basic difference of the  proof for the arbitrary geometry of $M$.
For the completeness of the proof we shall write the index $l$ in the further reasoning. This index should
be omitted for the two-point case.

The plan of the further proof is the following. First the basic reasoning will be given. After that we define the partitions of $\Lambda_n$ on the sets
$\Lambda_{nile}$ for the arbitrary geometry of $M$.  The basic reasoning is given on the set of events $A_{1n}$  such that
\begin{equation}\label{3.8b}
P(A_{1n}) = 1 + O(nb_n^{2+\lambda}).
\end{equation}
  The definition of the set $A_{1n}$ is rather cumbersome. To simplify the understanding
of the proof we have postponed the
definition of the set $A_{1n}$ to the end of section.

 For each $\theta \in \Lambda_{nile}$ denote
$$
S_{n\theta} = \sum_{s=1}^n \xi_{s}(\theta) - \theta'\sum_{s=1}^n
\tau_{s} + 2n\rho ^2(0,\theta)
$$
 and define the events
$$
B_{n\theta} = \left\{X_1,\ldots,X_n: S_{n\theta} > \epsilon_{1n}\right\}
$$
where $\epsilon_{1n} \to 0, \epsilon_{1n}^{-2}c_{1n}^{-3}nb_n^{2+\lambda} \to 0$ as $n \to \infty$.

Denote $B_{nile} = \cup_{\theta \in \Lambda_{nile}}B_{n\theta}$. For any  $\theta_{nij} \in \Theta_{nile}$
denote
$ B_{ni}(\theta_{nij})= \cup_{\theta\in\Lambda(\theta_{nij}) } B_{n\theta}$.

We have
\begin{equation}
\label{3.21}
\begin{split}&
\inf_{\hat\theta_n} \sum_{\theta \in \Lambda_{nile}}
P_\theta(\hat\theta_n - \theta \notin b_n \Omega, \psi_n \in n^{1/2}\Gamma_{ni})\\&\ge
\inf_{\hat\theta_n} \sum_{\theta \in \Lambda_{nile}}
E \left[\chi(\hat\theta_n - \theta \notin b_n
\Omega) \exp\left\{\sum_{s=1}^n
\xi_s(\theta)\right\}, \psi_n \in n^{1/2}\Gamma_{ni}, A_{1n} \right]\\& \ge E \left[\inf_t \sum_{\theta \in
\Lambda_{nile}} \chi(t - \theta \notin b_n \Omega)
\exp\left\{\theta \sum_{s=1}^n \tau_{s} -\frac {1}{2} n \theta'I \theta  + o(1)\right\},\right. \\&
\left.  \psi_n \in n^{1/2}\Gamma_{ni}, A_{nile} | A_{1n}\right]
P(A_{1n}) = R_n.
\end{split}
\end{equation}
Denote $\Delta_{n} = \exp\{\psi'_{n}\psi_{n}/2\}, y= y_\theta =
n^{1/2}\theta - \psi_{n}$. Then, using
 $nb_n\delta_n \to 0, nb_n^{2 +
\lambda} \to 0$ as $n \to \infty$, we get
\begin{equation}
\label{3.22}
\begin{split}&
(2l_n)^{-d}R_n \ge (2l_n)^{-d}E \left[\Delta_{n} \inf_t \sum_{\theta \in
\Lambda_{nile}} \chi(t - y_\theta - \psi_n \notin
n^{1/2}b_n \Omega)\exp\left\{-\frac{1}{2} y_\theta'I
y_\theta\right\},\right.\\& \left.\psi_n \in n^{1/2}\Gamma_{ni}, A_{nile}| A_{1n}\right] (1 +o(1)) \\&
= (2v_n)^{-d}E \left[\Delta_{n} \inf_t \int_{n^{1/2}K_{nile} -
\psi_n} \chi(t - y \notin n^{1/2}b_n
\Omega)\exp\left\{-\frac{1}{2} y'I y\right\}\,dy,\right.\\&\left.
 \psi_n \in n^{1/2}\Gamma_{ni}, A_{nile} | A_{1n}\}\right] (1 +o(1))\doteq
(2v_n)^{-d}I_{nile}(1 +o(1)).
\end{split}
\end{equation}
For each $\kappa \in (0,1)$  denote
$$
K_{ni\kappa}(\theta_{nij}) = \{x: x = \lambda x_{ni} +u +
\theta_{nij}, u=\{u_k\}_1^d, |u_k| \le ( c_{3n} - Cc_{2n})\delta_{1n},
$$
$$
u
\bot x_{ni},  \lambda \in R^1\} \cap K_{(1-\kappa)v_n},
$$
$$
K_{nile\kappa} = \cup_{\theta \in \Theta_{nile}} K_{ni\kappa}(\theta).
$$
If $\psi_n \in n^{1/2}\Gamma_{ni}\subset K_{\kappa v_n}$, then $n^{1/2} K_{nile\kappa}
\subset n^{1/2}K_{nile} - \psi_n$
and therefore
\begin{equation}
\label{3.23}
I_{nile} \ge U_{nile} \bar J_{nile} (1+o(1))
\end{equation}
with
$$
U_{nile}= E\left[\Delta_{n}, \psi_n \in \Gamma_{ni}, A_{nile}
|A_{1n}\right],
$$
$$
\bar J_{nile}\doteq \inf_t J_{nile}(t)\doteq  \inf_t
\int_{n^{1/2}K_{nile\kappa}} \chi(t - y \notin n^{1/2}b_n \Omega )
\exp\left\{-\frac{1}{2} y'I y \right\}dy.
$$
\begin{lemma}\label{l3.3}
\begin{equation}\label{3.24}
\bar J_{nile} = J_{nile}(0).
\end{equation}
\end{lemma}

Summing over $ l $ and $ e $, by (\ref {3.24}), we get
\begin{equation} \label{3.25}
\sum_{l=1}^{m_{3ni}}\sum_{e=1}^{m_{2nil}}\bar J_{nile\kappa}
\ge  P(I^{1/2}(\theta_0)\zeta \notin n^{1/2}b_n\Omega)(1+ o(1)).
\end{equation}
We have
\begin{equation}
\label{3.36}
\begin{split}&
U_{nile}= E\left[\Delta_{n}, \psi_n
\in n^{1/2}\Gamma_{ni} |A_{1n} \right]\\&-
E\left[\Delta_{n}, \psi_n \in
n^{1/2}\Gamma_{ni}, B_{nile} |A_{1n}\right] \doteq U_{1ni} - U_{2nile}.
\end{split}
\end{equation}
\begin{lemma}\label{l3.4} For all $i, 1 \le i \le m_n$
\begin{equation} \label{3.37}
U_{1ni}= \mbox{mes}(\Gamma_{ni})(1 + o(1)),
\end{equation}
\begin{equation}\label{3.38}
U_{2nile}= o(\mbox{mes}(\Gamma_{ni}))
\end{equation}
as $n \to \infty$.
\end{lemma}
Summing over $i$, by Lemma \ref{l3.4}, we get
\begin{equation}\label{3.39}
\sum_{i=1}^{m_n} U_{nile} \ge \mbox{mes}(K_{\kappa v_n})(1 +o(1))= (2\kappa v_n)^d(1 +o(1)).
\end{equation}
By (\ref{3.25},\ref{3.39}), we get
\begin{equation}\label{3.39a}
\sum_{i=1}^{m_n} \sum_{l=1}^{m_{3ni}}\sum_{e=1}^{m_{4ni}}
\bar J_{nile\kappa}U_{nile} \ge (2\kappa
v_n)^dP(I^{1/2}(\theta_0)\zeta \notin n^{1/2}b_n\Omega)(1 + o(1)).
\end{equation}
Since $\kappa, 0 < \kappa < 1$, is arbitrary, (\ref{3.8}), (\ref{3.21})-(\ref{3.23}),(\ref{3.39a}) together imply
Theorem \ref{t1}.

For the arbitrary geometry of the set $M$  the reasoning is the following.
Let us allocate in $M $ connectivity components $M_1, \ldots,
M_{s_1} $ having the greatest dimension. These components define
the asymptotic of  lower bound of risks. Denote
 $\tilde M = \cup_{i=1}^{s_1} M_i$.
Define the linear manifold $N$  having the smallest dimension $d_1$ such that $\tilde M \subset N$.
Define in $R^d$  the coordinate system, such that $N$  is induced the first  $d_1$ coordinates.
Denote $e_1,\ldots,e_d$ the vectors of the coordinate system.

Denote $y_{nij} \doteq y(\theta_{nij}) \doteq \{ x: x = \lambda x_{ni} + \theta_{nij}, \lambda >0\} \cap b_n\partial \Omega, 1 \le j \le m_{ni}$.
Define the sets $Y_{ni} = \{y : y = y_{nij}, 1\le j \le m_{1ni}\}$. We allocate in $Y _ {ni} $ the subset
$ \tilde Y _ {ni} $ of all points $y _ {nij} $ such that $K (\theta _ {nij}) \cap b_n\tilde M $ is not empty.

For each $y_{nij} \in \tilde Y_{ni}$ we set $z_{nij}\in  b_n\tilde M$ such that
\begin{equation}\label{3.5}
|y_{nij}- z_{nij}|=\inf_{z \in b_n\tilde M} |y_{nij} - z|.
\end{equation}
Define the set $\tilde Z_{ni}= \{z: z=z_{nij}, y_{nij} \in \tilde Y_{ni}\}$. Denote $ m_ {4ni} $ the number of points of
 $ \tilde Z_ {ni} $.

We split $\tilde Z_{ni}$ on subsets of points $\tilde Z_{nil} =
\{z_{nil1},\ldots,z_{nild_1}\}, 1\le l \le m_{3ni}$ such that the vectors
$ z_{nil1},\ldots,z_{nild_1}$ induce $N$. Note that $ t <d_1 $ points could  not enter in
  these partitions since $ m_ {4ni} $ may not be a multiple of
 $ d_1 $. However their exception is not essential for the further reasoning. Moreover, for the existence of
 such a partition we may have to define different constants $ c_{3n} $
 in the definition of different sets $ K_{nij} $. However, this does not affect significantly on the subsequent proof and
  we omit  the reasoning.

 For each $z_{nile}$ define the point $y_{nile},  y_{nile} \in \tilde Y_{ni}$ such that
$|y_{nile} - z_{nile}| \le c_{3n}\delta_{1n}$.

For each set $\tilde Z_{nil} \doteq \{z_{ni_1j_1},
\ldots, z_{ni_{d_1}j_{d_1}}\}=\{z_{nil1},\ldots,z_{nild_1}\}$
we make the following. For each point $\theta_{ni_sj_s},
1 \le s \le d_1$ we draw the linear manifold
$L_{i_sj_s}= \{z: z= \theta_{ni_sj_s}+ \lambda_1 e_{d_1+1} +
\ldots + \lambda_{d-d_1}e_d, \,\, \lambda_1,\ldots,\lambda_{d-d_1}
\in R^1\}$. We split $\Theta_{ni} \cap L_{i_sj_s}$ on the subsets
\begin{equation}\label{3.6a}
\begin{split}&
\Theta_{i_sj_s}(k_1,\ldots,k_{d-d_1})= \{\theta: \theta=
\theta_{ni_sj_s}  + (-1)^{t_1} 2k_1 c_{3n}\delta_{1n}e_{d_1+1}\\& +
\ldots + (-1)^{t_{d-d_1}} 2k_{d-d_1}c_{3n}\delta_{1n} e_d;\,\,
t_1,\ldots t_{d-d_1}= \pm 1\}
\end{split}
\end{equation}
where $0 \le k_1,\ldots,k_{d-d_1} < C_{1n}$ with $ C_{1n}c_{3n}c_{1n} \to \infty,
nb_n^{2+\lambda}C_{1n}^3c^{3}_{3n}c^{3}_{1n} \to 0$ as $n \to \infty$. Denote
$$
\tilde K_{i_sj_s}(k_1,\ldots,k_{d-d_1})= \cup_{\theta \in
\Theta_{i_sj_s}(k_1,\ldots,k_{d-d_1})} K(\theta).
$$
Denote $m_{2nil}(i_s,j_s)$ the number of sets $\tilde K_{i_sj_s}(k_1,\ldots,k_{d-d_1})$.

 Without loss of generality we can assume that
 $m_{2nil}(i_1,j_1)= m_{2nil}(i_2,j_2)=
 \ldots =m_{2nil}(i_d,j_d)\doteq m_{2nil}, 1 \le l \le m_{3ni}$.
This can always be achieved by making different constants
 $c_{3n}$ defining the sets $K_{nij}$.
Denote
\begin{equation}\label{3.6}
\bar K_{nil}(k_1,\ldots,k_{d-d_1}) = \cup_{s=1}^{d_1} \tilde
K_{i_sj_s}(k_1,\ldots,k_{d-d_1}).
\end{equation}
It will be convenient to number the sets $ \bar
K _ {nil} (k_1, \ldots, k _ {d-d_1}) $ denoting their $ \bar K _ {nil1},
\ldots,  \bar K_{nilm_{2nil}}$. Denote
\begin{equation}\label{3.7}
\Theta_{nile}=\Theta_{ni}\cap \bar K_{nile}, \quad \Lambda_{nile} = \bar
K_{nile}\cap \Lambda_n, \quad 1 \le e \le m_{2nil}.
\end{equation}
Thus $\Theta_{nile}$ contains $d_1 2^{d-d_1}$ points, that is,
$\Theta_{nile}= \{\theta_{sj}\}_{s=1,j=1}^{d-d_1,k}, k= 2^{d-d_1}$.

The further proof of Theorem \ref{t1} follows to the reasoning for the two-point $\{y,-y\}$ geometry of set $M$ given above.

Now the definition of the set $A_{1n} = A_{1nile}$ and the complementary set $B_{1n} = B_{1nile}=
D_{nile}\cup B_{4nile}\cup B_{3nile}$ will be given. The definitions of the sets
$D_{nile}, B_{4nile}, B_{3nile}$ are given bellow.

For all $s, 1 \le s \le n,$ denote $D_{ns}(\theta_{nij}) = \{X_s: f(X_s,0) \ne 0, f(X_s,\theta) = 0, \theta \ne 0,
\theta \in \Lambda_{n}(\theta_{nij}\}, D_n(\theta_{nij}) = \cup_{s=1}^n D_{ns}(\theta_{nij}),
D_{nile} = \cup_{\theta\in\Theta_{nile}} D_n(\theta)$.

Now we define the set $B_{2nile} \subset B_{4nile}$. For any $\theta_1,\theta_2 \in \Theta$ denote
$\eta_s(\theta_1,\theta_2) = g(X_s,\theta_1,\theta_2)$
with $1\le s \le n$.
Define the sets of  events $B_{2s}(\theta_1,\theta_2) = \{X_s:
|\eta_s(\theta_1,\theta_2)| \ge \epsilon\}, B_{2s}(\theta_2) = B_{2s}(0,\theta_2)$ with $0<\epsilon <\frac{1}{3}$.

For any $\theta \in \Theta_{nile}$ denote $B_{2nis}(\theta) = \cup_{\theta'
\in\Lambda_n(\theta)} B_{2s}(\theta')$,
 $B_{2ni}(\theta)= \cup_{s=1}^n B_{2nis}(\theta)$. Denote
$ B_{2niles} = \cup_{\theta \in\Theta_{nile}} B_{2nis}(\theta)$,
 $B_{2nile}= \cup_{s=1}^n  B_{2niles} $.

The estimates of $P(B_{2nile})$ are based on the "chaining method".
For simplicity we suppose that $l_n = 2^m$. This does not cause serious differences in the reasoning.
For each $\theta \in \Theta_{nile}$ we define the sets $\Psi_j= \Psi_j(\theta),
1\le j \le m$ of points $h_k = \theta +k\delta_{1n} e_1, h_k \in  \Lambda_{nile}, $ such that
 $ | k | $ is divisible by $ 2^{m-j} $ and is not divisible by $ 2 ^ {m-j +1},-l_{1n} \le k \le l_{1n} $.
Denote $\Psi_{m+1} = \Psi_{m+1}(\theta)= \Lambda_{n}(\theta)
\setminus \cup_{k=1}^m \Psi_k(\theta)$.   Denote
$\Psi_0(\theta) = \{\theta_0\}$. We say that the points $h \in
\Psi_j$ and $h_{1} \in \Psi_{j-1}$ are neighbors if $h_ {1}$ is
the nearest point of $ \Psi_ {j-1}$
 for $ h $.
For any $h \in \Psi_j$ we denote $\Pi(h)=\{  h_{1}: h_{1} \in
\Psi_{j-1}$ and $h, h_{1}-$ are neighbors $\}$.

For any $\theta \in \Theta_{nile}$ for each $h \in \Psi_j(\theta), 2 \le j \le m+1,$ and all $s, 1 \le s \le n$
define the events
$$
V_{hs}(\theta)= \{X_1: |\eta_s(h_{1},h)| > \epsilon j^{-2},
\eta_s(0,h_{1}) + 1 > \frac{1}{3} - \epsilon\sum_{k=0}^j k^{-2},  h_{1} \in \Pi(h)\}.
$$
Denote $$
B_{4nis}(\theta)=B_{2s}(\theta)\cup \cup_{2\le j\le
m+1}\cup_{h \in\Psi_j(\theta)}V_{hs}(\theta),\quad
B_{4niles}= \cup_{\theta \in \Theta_{nile}} B_{4nis}(\theta)
$$
and $
 B_{4nile}= \cup_{s=1}^n B_{4niles}(\theta).
$ It is clear that  $B_{2nis}(\theta)\subset B_{4nis}(\theta)$.
 \begin{lemma} \label{l3.1}
\begin{equation}\label{3.10}
P( B_{2nile}\cup D_{nile}) \le P( B_{4nile}\cup D_{nile})  = o(1).
\end{equation}
\end{lemma}

Define the event $B_{3ns} = \{X_s: |\tau_s| > \epsilon
v_n^{-1}\}$.
For any $\theta \in \Theta_{nile}$ for each $h \in \Psi_j(\theta), 1 \le j \le m+1,$ and all
$s, 1 \le s \le n$ define the events
$$
B_{3nhs} = \{X_s: |\tau_{sh} - \tau_s| > \epsilon b_n^{-1} 2^{j/2}\}.
$$
Denote
$$
B_{3nis}(\theta)=B_{3ns}\cup \cup_{2\le j\le
m+1}\cup_{h \in\Psi_j(\theta)}B_{3nhs}(\theta),\quad
B_{3niles}= \cup_{\theta \in \Theta_{nile}} B_{3nis}(\theta).
$$
and $
B_{3nile}(\theta) = \cup_{s=1}^n B_{3niles}
$
 \begin{lemma} \label{l3.2}
\begin{equation}\label{3.17}
P(B_{3nile} \cap A_{4nile}) = o(1).
\end{equation}
\end{lemma}
For any $\theta \in \Theta_{nile}$ denote $B_{1ns}(\theta) = B_{4ns}(\theta)\cup B_{3ns}(\theta) \cup D_{ns}(\theta)$.
Denote $B_{1n}(\theta) = \cup_{s=1}^n B_{1ns}(\theta), B_{1n}\doteq B_{1nile}=
\cup_{\theta\in\Theta_{nile}}B_{1n}(\theta).$

By Lemmas \ref{l3.1} and \ref{l3.2}, we get (\ref{3.8b}).

\section{Proofs of Lemmas \ref{l3.3} and \ref{l3.4}}

We begin with the proof of Lemma \ref{l3.4}.
The proof of (\ref{3.37}) is based on some version of Osypov-van Bahr Theorems  \cite{bahr,os}
on large deviation probabilities.

Let $Z$ be random vector in $R^d$ such that $E[Z] = 0, \mbox{Var}
(Z) = I$, where $I$  is unit matrix. Let $P(|Z| <\epsilon b_n^{-1}) = 1$, where $\epsilon > 0$ as $n
\to \infty$.  Suppose $E |Z|^{2+\lambda} <C < \infty$.
 Let $Z_1,\ldots,Z_n$ be independent copies of $Z$. Denote $S_n = n^{-1/2}(Z_1+\ldots+ Z_n)$.

Denote $\mu_n$ the probability measure of Gaussian random vector $\zeta$ with $E[\zeta] = 0$ and covariance matrix
 $nI$. For any Borel set $W$ denote $W_\delta$
$\delta$- vicinity of  $W, \delta>0$.

\begin{theorem}\label{t2} Let the set $W$ belong to a ball in $R^d$ having the radius $r =
o(\epsilon_n n^{1/2} b_n)$ where $\epsilon_n \to 0$ as $n \to 0$. Let $nb_n^2 \to \infty, nb_n^{2+\lambda} \to 0$ as $n \to \infty$.  Let $W = W_1 \setminus W_2$ where $W_1,W_2$ are
the convex sets. Then
\begin{equation}\label{3.40}
P(S_n \in W) = \mu_n(W)(1 + O(b_n^\lambda)) + O(b_n^\lambda)\mu_n(W_{c_{n}})
\end{equation}
where $c_{n} =  o(n^{-1/2}b_n^{\lambda-1})$.
\end{theorem}

  The differences in the statements of Theorem \ref{t2} and Osypov - van Bahr Theorem \cite{bahr,os} are caused
the differences in the assumptions.
In  \cite{bahr,os} the results have been proved if $ E[\exp\{c|Z|\}] < \infty$.

Let us check up that the assumptions of Theorem  \ref {t2}  are
fulfilled for the random vector $Z =
I^{-1/2}(\theta_0)\tau\chi(A_{1n1})$.
\begin{lemma}\label{l3.5}
\begin{equation}\label{3.41}
E[\tau, A_{1n1}] = O(b_n^{1+\lambda}),
\end{equation}
\begin{equation}\label{3.42}
E[\tau\tau', A_{1n1}] = I(\theta_{0})+
O(b_n^{\lambda}).
\end{equation}
\end{lemma}
Lemma \ref{l3.5} and Theorem \ref{t2} imply (\ref{3.37}).

\begin{lemma} \label{l3.6} Uniformly in $\theta \in \Lambda_{nile}$
\begin{equation}\label{3.20}
E_\theta [S_{n\theta}| A_{1n}] = o(1).
\end{equation}
\end{lemma}
Let $\epsilon_{1n}$  be such that
\begin{equation}\label{3.20aa}
\sup_{\theta\in\Lambda_{nile}} |E[S_{n\theta}|A_{1n}] \le \frac{\epsilon_{1n}}{4}.
\end{equation}
Let $h  \in \Psi_j,  h_1
\in \Pi(h), 2 \le j \le m+1$. We have
\begin{equation}\label{3.44}
S_{nh} - E[S_{nh}|A_{1n}]  =
S_{nh_{1}}+S_{1nh}+ S_{2nh} - E[S_{nh_{1}}+S_{1nh}+
S_{2nh}|A_{1n}]
\end{equation}
where
\begin{equation}\label{3.45}
S_{1nh}=\sum_{s=1}^n\xi_s(h_{1},h) - \bar h'\sum_{s=1}^n
\tau_{sh_{1}} ,
\end{equation}
\begin{equation}\label{3.46}
S_{2nh} = \bar h'\sum_{s=1}^n (\tau_{sh_{1}} - \tau_s
)
\end{equation}
with $\bar h = h - h_1$.

Denote
$$
B_{0n}= \{X_1,\ldots,X_n: \sup_{h \in \Psi_1} S_{nh} >
\epsilon_{1n}/4\}.
$$
For any $h \in \Psi_j, 2\le j \le m+1$ denote
$$
B_{5nh} = \{X_1,\ldots,X_n: j^2 (S_{1nh}- E[ S_{1nh}|A_{1n}])
> \epsilon_{1n}/4\},
$$
$$
B_{6nh} = \{X_1,\ldots,X_n: j^2 (S_{2nh} - E [S_{2nh}
|A_{1n}])
> \epsilon_{1n}/4\}.
$$
Denote $ B_{n} = B_{0n}\cup (\cup_{\theta\in \Lambda_{nile}
\backslash \Psi_1} (B_{5n\theta}\cup B_{6n\theta}))$.
Note that  $ B_{n}\supseteq B_{nile}$. Hence
\begin{equation}\label{3.47}
U_{2nile} \le U_{3nile}\doteq E\left[\Delta_{n}, \psi_n \in
n^{1/2}\Gamma_{ni},  B_{n} |A_{1n}\right].
\end{equation}

Denote  $r_{ni}= \inf_{x\in \Gamma_{ni}} |x|$.
We have
\begin{equation}\label{3.48}
U_{3nile} \le C \exp\{n r_{ni}^{2}/2\} \left(V_{0n} +
\sum_{\theta\in \Lambda_{1nile} }(V_{5n\theta }
+V_{6n\theta})\right)
\end{equation}
where $\Lambda_{1nile} = \Lambda_{nile}\setminus \Theta_{nile}$,
\begin{equation}\label{3.49}
V_{en\theta } = P \left(\psi_n \in
n^{1/2}\Gamma_{ni}, B_{en\theta} \,|\,A_{1n} \right), \quad e= 5,6,
\end{equation}
\begin{equation}\label{3.50}
V_{0n} = P  \left(\psi_n \in
n^{1/2}\Gamma_{ni}, B_{0n} \,|\,A_{1n}
\right).
\end{equation}
\begin{lemma} \label{l3.7}  Let $\zeta$  Gaussian random vector having the covariance matrix
 $I(\theta_0)$ and let $E[\zeta] =  0$. Then for any  $h \in
\Psi_j, h_{1} \in \Pi(h)$
\begin{equation}\label{3.51}
V_{0n} \le C n b_n^{2+\lambda} \epsilon_{1n}^{-2}P(\zeta \in n^{1/2}\Gamma_{ni}),
\end{equation}
\begin{equation}\label{3.52}
V_{5nh} \le C n |\bar h|^{2} b_n^\lambda \epsilon_{1n}^{-2} j^4 P(\zeta \in
n^{1/2}\Gamma_{ni}),
\end{equation}
\begin{equation}\label{3.53}
V_{6nh} \le C n|\bar h|^{2}  b_n^\lambda \epsilon_{1n}^{-2} j^4 P(\zeta \in
n^{1/2}\Gamma_{ni}).
\end{equation}
\end{lemma}
The number of points $\Psi_j, 1 \le j \le m,$ equals $2^j$ and, if $h \in \Psi_j$, then $\bar h = b_n2^{-j}$.
The number of points $\Psi_{m+1}$ equals $C c_{3n}^{d-1}2^m$ and, if $h \in \Psi_{m+1}$,
then $|\bar h| \le Cc_{3n}\delta_{1n}$.
Hence, by Lemma \ref{l3.7},  we get
\begin{equation}\label{3.54a}
\begin{split}&
U_{3nile} \le Cn\epsilon_{1n}^{-2} \exp\{nr_{ni}^2/2\} P(\zeta \in
n^{1/2}\Gamma_{ni})\\&\times \left(b_n^{2+\lambda} +
b_n^\lambda\left(\sum_{j=1}^m 2^j (b_n2^{-j})^2 j^4 + c_{3n}^{d+1} m^4 2^m
\delta_{1n}^2\right)\right).
\end{split}
\end{equation}
Note that $m$ satisfies $\delta_{1n} = v_n 2^{-m}$ or
 $2^m = C c_{1n}^{-1}n b_n^2(1+ o(1))$. Hence
\begin{equation}\label{3.54b}
n\epsilon_{1n}^{-2} b_n^\lambda c_{3n}^{d+1}  m^4 2^m \delta_{1n}^2 = C n\epsilon_{1n}^{-2}
b_n^\lambda c_{3n}^{d+1}  c_{1n}^{-1} nb_n^2 m^4 c_{1n}^{-2}
n^{-2} b_n^{-2} = C \epsilon_{1n}^{-2} b_n^\lambda c_{3n}^{d+1}  c_{1n}^{-3} m^4 = o(1).
\end{equation}
By (\ref{3.54a}, \ref{3.54b}), we get
\begin{equation}\label{3.54}
U_{3nile} = o(mes(\Gamma_{ni})).
\end{equation}
By (\ref{3.47}) and (\ref{3.54}), we get (\ref{3.38}).

Proof of Lemma \ref{l3.7} is based on Theorem \ref{t3}.

\begin{theorem} \label{t3} Let we be given a random vector $V = (X,Z)$ where random variable $X$  and random vector
$Z = (Z_1,\ldots,Z_d)$ are
such that $E[V] = 0$. Let
\begin{equation}\label{3.55}
P(|X| < \epsilon) = 1, \quad  E[|X|^2]  < Cb_n^{2+\lambda},
\end{equation}
\begin{equation}\label{3.56}
P(|Z| < \epsilon b_n^{-1}) = 1,
 \quad E [|Z|^{2+\lambda}] < C < \infty,
\end{equation}
\begin{equation}\label{3.57}
E [XZ_k] = O(b_n^{1 +\lambda}), \quad 1\le k \le d
\end{equation}
with $0<\epsilon<1$. Suppose the covariance matrix of random vector $Z$ is positively definite.

Let  $V_1 = (X_1,Z_1), \ldots, V_n = (X_n,Z_n)$  be independent copies of random vector
 $V$. Let $U $ be a bounded set in
$R^d $ being a difference of two convex sets.

 Denote $S_{nX} =  n^{-1/2}(X_1 +
\ldots + X_n)$ and $ S_n
= n^{-1/2}(Z_1 + \ldots + Z_n)$. Denote $Y$ the Gaussian random vector having the same covariance matrix as the random vector $Z$.

Then, for the sufficiently large  $n$,
\begin{equation}\label{3.58}
I \doteq P(S_{nX} > \epsilon_{1n}, S_n \in nb_n v + r_n  U) \le
CP(S_{nX} > \epsilon_{1n}) P(Y \in nb_n v + r_n  U)
\end{equation}
where $\epsilon_{1n}, r_n$ are chosen so that
 $nb_n^{2+\lambda}c_{n1}^{-3}\epsilon_{1n}^{-2}\to 0$ as
$n \to \infty$ and $r_n > c_{n1} n^{-1/2}b_n^{-1}$.
\end{theorem}
It is clear  that $\epsilon_{1n}, r_n$ can be chosen such that
$\epsilon_{1n} \to 0, r_n n^{1/2}b_n \to 0$ as $n \to \infty$.
 In the proof of (\ref{3.52},\ref{3.53})
we suppose that $\epsilon_{1n}$ and $r_n$ satisfy these assumptions.

For  the estimates of  $V_{5nh}$ in (\ref{3.52}) we implement Theorem \ref{t3} with $Z=\tau$ and
$$
X=\varphi(h_1,h) =\xi(h_1,h) - \bar h' \tau_{h_1} - \sum_{k=1}^d\rho_{kh_1h}\tau_k.
$$
Here $\tau =\{\tau_k\}_{k=1}^d$ and $\rho_{h_1h}= \{\rho_{kh_1h}\}_{k=1}^d=r_{h_1h}(E[\tau\tau'|A_{1n1}])^{-1}$ with $r_{h_1h} = \{r_{kh_1h}\}_{k=1}^d,
 r_{kh_1h} = E[(\xi(h_1,h) -
\bar h' \tau_{h_1})\tau_k| A_{1n1}].$

Thus  $S_{1nh}$ is replaced with
$$
S_{nx} = S_{1nh} -  \sum_{s=1}^n \sum_{k=1}^d\rho_{kh_1h}\tau_{ks} = \sum_{s=1}^n \varphi_s(h_1,h).
$$

It is easy to see that  $E[\varphi(h_1,h)\tau_k| A_{1n1}]= 0, 1\le k \le d$. This implies  (\ref{3.57}).

Now we show that
\begin{equation}\label{3.58a}
\sum_{s=1}^n\sum_{k=1}^d \rho_{kh_1h}\tau_{ks} = o(1)
\end{equation}
if $\psi_n \in n^{1/2}\Gamma_{ni}$
This justifies such a replacement.

By Lemma \ref{l3.10} given bellow,
$|r_{kh_1h}| \le C |\bar h|^{1+\lambda/2}$, if $2\le k \le d$. Hence, since $\psi_n \in n^{1/2} \Gamma_{ni}$,
\begin{equation}\label{3.59}
r_{kh_1h} \sum_{s=1}^n \tau_{ks}=O( |\bar h|^{1+\lambda/2}
b_n^{-1}) = o(1)
\end{equation}
with $2\le k \le d$.

\begin{lemma}\label{l3.10}
Let  $h \in \Psi_j(\theta), 1\le j \le m+1$, $h_1 \in \Pi(h)$ and let $v \perp \bar h, u \in R^d$. Then
\begin{equation} \label{3.64}
E[(\xi(h_1,h) - \bar h' \tau_{h_1}) (v'\tau), A_{1n1}] = O(|v| |\bar h |^{1+\lambda/2}).
\end{equation}
\end{lemma}
By Lemma \ref{l3.11} given bellow $|r_{1h_1h}| \le C |\bar h|b_n^{\lambda}$. Hence, since $\psi_n \in n^{1/2} \Gamma_{ni}$,
\begin{equation}\label{3.60}
r_{1h_1h} \sum_{s=1}^n \tau_{1s}=O( n|\bar h|
b_n^{1+\lambda}) = o(1).
\end{equation}
By (\ref{2.4}), (\ref{3.59}), (\ref{3.60}), we get (\ref{3.58a}).
\begin{lemma}\label{l3.11}
Let  $h \in \Psi_j(\theta), 1\le j \le m+1$, $h_1 \in \Pi(h)$ and let
 $v \parallel \bar h$. Then
\begin{equation} \label{3.65}
E[(\xi(h_1,h) - \bar h'\tau_{h_1}) (v'\tau), A_{1n1}] =
O(|v| |\bar h | b_n^{\lambda}).
\end{equation}
\end{lemma}
Note that
\begin{equation}\label{3.63u}
2\eta(h_1,h) - 2\eta^2(h_1,h) \le \xi(h_1,h) \le 2 \eta(h_1,h) < 2\epsilon
\end{equation}
if $A_{1n1}$ holds.

By (\ref{3.63u}) and Lemma \ref{l3.8} given bellow, we get (\ref{3.55}).
\begin{lemma}\label{l3.8} For all $\theta \in \Lambda_{nile}$
\begin{equation}\label{3.61}
E[(\xi(\theta) - \theta'\tau)^2, A_{1n1}] =
O(|\theta|^{2+\lambda}).
\end{equation}
Let $h \in \Psi_j(\theta), 1\le j \le m+1$ и $h_1 \in \Pi(h)$. Then
\begin{equation}\label{3.62}
E[(\xi(h_1,h) - \bar h' \tau_{h_1})^2, A_{1n1}] = O(|\bar h |^{2+\lambda}).
\end{equation}
\end{lemma}
This completes the proof of (\ref{3.52}).

The proof of
 (\ref{3.51}) is akin to the proof of (\ref{3.52}) and is omitted.

For  the estimates of  $V_{6nh}$ in (\ref{3.53}) we choose $Z=\tau$ and
$$
X\doteq
\bar h'(\tau_{h_1} - \tau)- \sum_{k=1}^d
\bar\rho_{kh_1h}\tau_k.
$$
Here $\tau =\{\tau_k\}_{k=1}^d$ and $\bar\rho_{kh_1h}= \{\bar\rho_{kh_1h}\}_{k=1}^d=
\bar r_{h_1h}(E[\tau\tau'|A_{1n1}])^{-1}$ with $\bar r_{h_1h} = \{\bar r_{kh_1h}\}_{k=1}^d,
\bar r_{kh_1h} = E[\bar h'(\tau_{h_1}- \tau)\tau_k| A_{1n1}], 1 \le k \le d
$.

Using the same reasoning as in the proof of (\ref{3.52}) and Lemmas \ref{l3.9}, \ref{l3.12} given bellow we get (\ref{3.53}).
\begin{lemma} \label{l3.9} Let $u,h \in R^d$.  Then
\begin{equation}\label{3.63}
E[(u'(\tau - \tau_h))^2, A_{1n1}] = O(|u|^2 |h|^\lambda).
\end{equation}
\end{lemma}
\begin{lemma} \label{l3.12} Let $h \in \Psi_j(\theta), 1 \le j \le m+1, h_1 \in \Pi(h).$ Let $v \perp \bar h, v \in R^d$. Then
\begin{equation}\label{3.63c}
E[\bar h'(\tau_{h_1} - \tau)(v'\tau), A_{1n1}] = O(|v| |\bar h| |h_1|^{\lambda/2}).
\end{equation}
If $v \parallel \bar h$,
\begin{equation}\label{3.63d}
E[\bar h'(\tau_{h_1} - \tau)(v'\tau), A_{1n1}] = O(|v| |\bar h| |h_1|^{\lambda}).
\end{equation}
\end{lemma}
{\it Proof of Lemma \ref{l3.3}}. The set $ \Lambda_{nile} $ is defined by the set of the points
 $\Theta_{nile} = \{\theta_{sj}\}_{s,j=1}^{d_1,k}, k = 2^{d-d_1}$.
The reasoning first will be given for  $|t| < c < \infty$.
Denote $n^{1/2}y_{sj}(t) \in (n^{1/2}b_n\partial \Omega - t)\cap (n^{1/2} K(\theta_{sj}))$ the point
in which $n ^ {1/2} y_{sj} =n ^ {1/2} y (\theta _{sj}) $ will pass  at the shift $t $.
Denote  $n^{1/2}y_{s+d_1,j}(t) \in (n^{1/2}b_n\partial \Omega - t)\cap (n^{1/2} K(\theta_{sj}))$ the point
in which $n^{1/2}y_{d_1+s}=-n^{1/2}y_{sj}$ will pass  at the shift $t$.
\begin{lemma}\label{l4.21} There holds
\begin{equation}\label{3.26}
\sum_{s=1}^{2d_1} \sum_{j=1}^k \exp\left\{-\frac{1}{2}n|y_{sj}(t)|^2\right\}
\ge 2\sum_{s=1}^{d_1} \sum_{j=1}^k \exp\left\{-\frac{1}{2}n|y_{sj}|^2\right\}.
\end{equation}
\end{lemma}
{\it Proof of Lemma \ref{l4.21}}. For a while we fix $ s \le d_1 $ and $ j $.
 We slightly modify the coordinate system for the further reasoning.
 Suppose $x_{ni} =
(1,\beta_2,\ldots,\beta_d)$ and  $y_{sj}=
(b_n,0,\ldots,0,\delta_{d_{1}+1,n}n^{-1/2}$,  $\ldots,
\delta_{dn}n^{-1/2})(1 + o(n^{-1/2}b_n^{-1}))$ with $\delta_{kn} \in R^1, d_1+1 \le k \le d$.

Define the line $y = n^{1/2}(y_{sj} + u  x_{ni}), \quad u \in R^1$, that is,
$$
y_1 = n^{1/2}b_n + u, y_2 =  \beta_2u,\ldots,x_{d_1} =
\beta_{d_1}u,
$$
$$
y_{d_1+1}= \delta_{d_1+1,n} + \beta_{d_1+1} u, \ldots,y_{d}=
\delta_{d,n} + \beta_{d} u,
 \quad |\delta_{kn}| < C,  d_1+1\le k \le d,  u \in R^1.
$$
Denote  $\delta_{kn} = 0$ for $1 < k \le d_1$.

Since the reasoning is given in a sufficiently small
 vicinity of point
 $n^{1/2}y_{sj}$
 the surface $n^{1/2}b_n\partial\Omega$
 admits  the approximation in this vicinity by an ellipsoid
$$
(x_1 - n^{1/2}b_n)^2 + \alpha_2 x_2^2+\ldots + \alpha_d x_d^2 =  nb_n^2
$$
where $-\alpha_2,\ldots,-\alpha_d$ are the principal curvatures of the surface
$\partial \Omega$ at the point $(1,0,\ldots,0)$. Thus, in  the further reasoning,
we can replace the set $n^{1/2}b_n\partial\Omega$
with the ellipsoid.
After the shift $ t = (t_1, \ldots, t_d) $ the ellipsoid is defined by the equation
$$
(x_1 - n^{1/2}b_n+ t_1)^2 + \alpha_2 (x_2+t_2)^2+\ldots + \alpha_d (x_d + t_d)^2 =  nb_n^2
$$
 and intersects the line $ y = n^{1 / 2} (\theta_ {sj} + ux_ {ni}), u \in R ^ 1 $
at the point $n^{1/2}y_{sj}(t)$ having the coordinates
\begin{equation}\label{3.27}
n^{1/2}y_1(t) = n^{1/2}b_n-  t_1+ \omega_{1n},  n^{1/2}y_k(t) = \delta_{kn}- \beta_2 t_1 +\beta_2\omega_{1n},\quad
1<k \le d.
\end{equation}
with
\begin{equation}\label{3.28}
\omega_{1n} = - (2n^{1/2}b_n)^{-1}(\alpha_2(\delta_{2n} + t_2-\beta_2 t_1)^2 + \ldots +
\alpha_d(\delta_{dn} + t_d - \beta_d t_1)^2)(1+o(1)).
\end{equation}
Arguing similarly we get that the ellipsoid intersects the line
$y =n^{1/2}( -y_{sj} + u x_{ni}),
u \in R^1$
at the point $n^{1/2}y_{s+d_1,j}(t)$ having the coordinates
\begin{equation}\label{3.29}
n^{1/2}y'_1(t) = -n^{1/2}b_n-  t_1+ \omega_{2n},\quad  n^{1/2}y'_s(t) =- \delta_{kn} - \beta_k t_1 +\beta_k\omega_{2n}
\quad 1 < k \le d_1
\end{equation}
with
\begin{equation}\label{3.30}
\omega_{2n} = (2n^{1/2}b_n)^{-1}(\alpha_2(-\delta_{2n} + t_2-\beta_2 t_1)^2 +
\ldots + \alpha_d(-\delta_{dn} + t_d - \beta_d t_1)^2)(1+o(1)).
\end{equation}
Substituting (\ref{3.27}, \ref{3.29}) in  (\ref{3.26}) we find that, if $t_1 >>>  n^{-1/2}b_n^{-1}$, then
 $$
\max\{\exp\{-n(y_1(t)^2)/2\}, \exp\{-n(y'_1(t)^2)/2\}\} >>>
\exp\{-(nb_n^2+\delta_{d_1+1}^2\ldots+\delta_d^2)/2\}.
$$
Thus we can suppose   $t_1 < cn^{-1/2}b_n^{-1}$   and
neglect the addendums $\beta_it_1, 2\le i \le d$ in (\ref{3.28},\ref{3.30}).

Using (\ref{3.27}, \ref{3.29})
, we get
\begin{equation}\label{3.31}
\begin{split}&
\exp\left\{-\frac{1}{2}n|y_{sj}(t)|^2\right\}+ \exp\left\{-\frac{1}{2}n|y_{s+d_1,j}(t)|^2\right\}\\&=
\exp\{-n|y_{sj}|^2/2\}\left(\exp\left\{n^{1/2}b_nt_1 +
\sum_{k=d_1+1}^{d}\alpha_k t_k\delta_{kn}\right\} \right.\\&  \left.
+\exp\left\{-n^{1/2}b_nt_1- \sum_{k=d_1+1}^{d}\alpha_k t_k\delta_{kn}\right\}\right)
\exp\left\{\frac{1}{2}\sum_{k=d_1+1}^{d}\alpha_kt_k^2\right\}
(1+ o(1)).
\end{split}
\end{equation}
Taking the points $ y_ {sj}, 1 \le j \le 2^{d-d_1},$
 with all possible values
 $\pm \delta_{kn},d_1 < k \le d$ and
 summing up for them
 $\exp\{-\frac{|y^2_{sj}(t)|^2}{2}\}$  we get
\begin{equation}\label{3.32}
\begin{split}&
\exp\left\{-\frac{nb_n^2 + \delta_{d_1+1,n}^2 + \ldots \delta_{dn}^2}{2}\right\}\\&\times
(\exp\{n^{1/2}b_nt_1\} + \exp\{-n^{1/2}b_nt_1\})\\& \times
\prod_{k=d_1+1}^{d} (\exp\{\alpha_k t_k\delta_{kn}\} +
\exp\{-\alpha_k t_k\delta_{kn}\} )(1 + o(1)).
\end{split}
\end{equation}
Since $\exp\{v\} + \exp\{-v\} -2 \ge 0$ with $v \in R^1$, then
(\ref{3.32}) implies (\ref{3.26}) for $|t | < C$.

In essence, we have considered only the case
 $u = 0$. Any point $y_u = n^{1/2}(y_{sj} + ux_{ni}), 0 < u <<< 1$, pass after the shift $t$ at the point
$n^{1/2}(y_{sj}(t) + ux_{ni}) \in (R^d\setminus (n^{1/2}b_n \Omega - t))\cap (n^{1/2} K(\theta_{sj}))$. Thus for any point
$y_u, 0 < u <<< 1$ we can write a similar inequality  (\ref{3.26}).
Since the shift $t$ is negligible,
\begin{equation}\label{3.33}
\mbox{mes}((n^{1/2}b_n\partial\Omega) \cap K(\theta_{sj}) )=
\mbox{mes}((n^{1/2}b_n\partial\Omega- t) \cap K(\theta_{sj}))(1 + o(1)).
\end{equation}
This implies $\bar J_{nile}(t) \ge J_{nile}(0)$.

Let us consider the case $ c <<  |t| << Cn^{1/2}b_n$. Note that, since all the principal curvatures in all points of $\partial \Omega$
 are negative, we can conclude $ n^{1/2} b_n \Omega $ into an ellipsoid
$$
 \Xi=\{x=\{x_i\}_{i=1}^d :  x_1^2+ \ldots + x_{d_1}^2 +
 \bar\alpha_{d_1+1} x_{d_1+1}^2 +\ldots+\bar\alpha_{d}x_d^2=
nb_n^2\}
$$
passing through the points $ y_{nile}$ and $ -y_{nile}, 1\le e \le d_1$ and such that
 $\bar\alpha_k < 1, d_1+1 \le k \le d$.
Denote $ y_{sj} (t) \in (n^{1/2} b_n \partial \Omega - t) \cap
\{y: y = \theta_{sj} +
  x_{ni} u, u \in R^1 \} $ and $ \bar y_{sj} (t) \in
 (\Xi-t) \cap \{ y: y = \theta_ {sj} + x_ {ni} u, u \in R ^ 1 \} $
the point in which the $ y_{sj}$ will pass at the shift $ t $.

It is easy to see
\begin{equation}\label{3.34}
\sum_{s=1}^{2d_1} \sum_{j=1}^k\exp\left\{-\frac{|y_{sj}(t)|^2}{2}\right\}\ge
\sum_{s=1}^{2d_1} \sum_{j=1}^k\exp\left\{-\frac{|\bar y_{sj}(t)|^2}{2}\right\}.
\end{equation}
For the points  $\bar y_{sj}(t)$ we can make estimates similar to the case $|t| < C<\infty$ and  can get
\begin{equation}\label{3.35}
\sum_{s=1}^{2d_1} \sum_{j=1}^k\exp\left\{-\frac{|\bar y_{nils}(t)|^2}{2}\right\}\ge
\sum_{s=1}^{2d_1} \sum_{j=1}^k\exp\left\{-\frac{| y_{nils}|^2}{2}\right\}.
\end{equation}
The statement  (\ref{3.35}) implies $J(t) > J(0)$ for  $ c << |t|
<<Cn^{1/2}b_n$.

Finally, after the shift $ t, |t| \asymp n ^ {1 / 2} b_n $ one of the points $ y_{nile}$ or $ -y_{nile}, 1\le e \le d_1$
will be located at a distance having the
 order $ n^{1/2} b_n $ outside  the ellipsoid $\Xi$ and hence outside  $ n^{1/2}b_n\Omega $.
 This implies $ J (t)> J (0) $.
 \section{Proofs of Theorems \ref{t2} and \ref{t3} }
The proof of Theorem \ref{t2} contains only some different technical details in comparison with the proof of similar Theorem
in \cite{os}.
The proof of Theorem \ref{t3} is based on a fairly new analytical technique
(see  \cite{be, er01}) and is more interesting.
Thus we begin with the proof of Theorem \ref{t3}.

{\it Proof of Theorem \ref{t3}}. We begin with auxillary estimates of moments of random variable  $X$ and random vector $Z$. We have
\begin{equation}\label{4.12}
E [|X| |Z|^2] \le (E |X|^{\frac{2+\lambda}{\lambda}})^{\frac{\lambda}{2+\lambda}}
(E|Z|^{2+\lambda})^{\frac{2}{2+\lambda}} \le
C (E[ X^2])^{\frac{\lambda}{2+\lambda}} \le C b_n^\lambda,
\end{equation}
\begin{equation}\label{4.13}
E[X^2 |Z|] \le  C b_n^{-1} E[ X^2] \le C b_n^{1+\lambda},
\end{equation}
\begin{equation}\label{4.14}
E[X^2 |Z|^2] \le  C b_n^{-2} E[ X^2] \le C b_n^{\lambda},
\end{equation}
\begin{equation}\label{4.15}
E[ X^2 |Z|^3] \le  C b_n^{-3} E[ X^2] \le C b_n^{\lambda-1},
\end{equation}
\begin{equation}\label{4.16}
E[ X^2 |Z|^3] \le  C  E[ |Z|^3] \le C b_n^{\lambda-1} E[ |Z|^{2+\lambda}] \le  C b_n^{\lambda-1}.
\end{equation}
 For each $x = \{x_1,\ldots,x_d\} \in R^d$ denote
$||x|| = \max_{1\le i \le d} |x_i|$. For any $z \in R^d$ and any  $A \subset R^d$
denote $||A-z||= \inf_{x\in A} ||x-z||$. For any $\epsilon > 0$ denote
$A_\epsilon = \{x : ||A-x|| \le \epsilon, x \in R^d\}$.

Define twice continuously differential functions $f_{1n}: R^1 \to R^1$
such that
$$
f_{1n}(x)=
  \begin{cases}
1 \quad \mbox{if}  \quad |x| > \epsilon_{1n}
\\
0 \quad \mbox{if}  \quad |x| < \epsilon_{1n}/2
  \end{cases}
$$
and $0 \le f_{1n}(x) \le 1, \left|\frac{\partial f_{1n}(x)}{\partial x_{i_1}\partial x_{i_2}}\right| \le
C\epsilon_{1n}^{-2},
1\le i_1,i_2 \le d, x \in R^d$.

Denote $c_{n}= c_{n1} n^{-1/2}b_n^{-1}$. We  slightly modify the setup of Theorem \ref{t3} in the proof.
The reasoning will be given with
 $r_n=1$.  Theorem \ref{t3} follows from the reasoning if we put $r_n =c_{n}$.

Define three- times continuously differentiable functions
$f_{2n} : R^d \to R^1$ such that
$$
f_{2n}(x)=
  \begin{cases}
1 \quad \mbox{if} \quad x \in n^{1/2}b_n v +  U
\\
0 \quad \mbox{if} \quad  x \notin n^{1/2}b_n v +
U_{c_{n}}
  \end{cases}
$$
and $0 \le f_{2n}(x) \le 1, |\frac{\partial^3 f_{2n}(x)}{\partial x_{i_1}
\partial x_{i_2}\partial x_{i_3}}| \le Cc_{n}^{-3}, 1 \le i_1,i_2,i_3 \le d$ if
 $ x \in R^d$.

Denote
$$
S_{knX} =X_1+\ldots+ X_{k-1}+ X_{k+1}+\ldots+X_n,
$$
$$
W_{kn} = n^{-1/2}(Z_1+\ldots + Z_{k-1}+Y_{k+1}+\ldots+Y_n).
$$
Hereafter $ Y_1,\ldots, Y_n$ are independent copies of random vector $Y$.
 Random variables $Y, Y_1,\ldots, Y_n$ do not depend on $X_1,\ldots, X_n, Z_1,\ldots,Z_n$.

For any $\gamma > 0$ denote
$$
G_n(\gamma) = sup \,\,E[ f_{1n}(S_{nX}), S_{nZ} \in n^{1/2}b_nv +  U_\gamma]
$$
where the supremum is taken over all distributions of  $(X,Z)$ satisfying the assumptions of Theorem \ref{t3}.

\begin{lemma} \label{l4.2} Let assumptions of Theorem \ref{t3} be satisfied. Then
\begin{equation}\label{4.17}
\begin{split}&
E[f_{1n}(S_{nX}), S_{nZ} \in n^{1/2}b_n v + U] \\& \le
E[f_{1n}(S_{nX})] P(Y \in n^{1/2}b_nv + U_{c_n}) +
Cnb_n^{2+\lambda} c_{n1}^{-3}\epsilon_{1n}^{-2} G_{n-1}(\gamma_n)
\end{split}
\end{equation}
for $n > n_0$.
Here $\gamma_n = \epsilon b_n^{-1}(n-1)^{-1/2}  + (n(n-1)^{-1/2}b_n - (n-1)^{1/2}b_{n-1}) + C/n + c_n$ where $C$ depends on $U$.
\end{lemma}
{\it Proof of Lemma \ref{l4.2}.} We have
\begin{equation}\label{4.18}
E[f_{1n}(S_{nX})f_{2n}(S_{nZ})] \le E[f_{1n}(S_{nX})f_{2n}(Y)] + \Delta
\end{equation}
where
\begin{equation}\label{4.19}
\Delta = |E[f_{1n}(S_{nX})f_{2n}(S_{nZ})] - E[f_{1n}(S_{nX})f_{2n}(Y)]|.
\end{equation}
It is clear that $\Delta \le \Delta_1+ \ldots + \Delta_n$ where
\begin{equation}\label{4.20}
\Delta_k = |E[f_{1n}(S_{knX}+ X_k)f_{2n}(W_{kn} + n^{-1/2}Z_k)] -
E[f_{1n}(S_{knX}+ X_k)f_{2n}(W_{kn} + n^{-1/2}Y)]|
\end{equation}
for $1 \le k \le n$.

Expanding $f_{1n}$ and $f_{2n}$ in the Taylor series, we get
\begin{equation}\label{4.21}
\begin{split}&
\Delta_k = |E[f_{1n}(S_{knX}+ X_k)(f_{2n}(W_{kn} + n^{-1/2}Z)] -
f_{2n}(W_{kn} + n^{-1/2}Y))]|\\&
\le \left|E\left[\left(f_{1n}(S_{knX}) + f'_{1n}(S_{knX}) X_k +
\frac{1}{2}\int_0^1 f''_{1n}(S_{knX} + \omega X_k) (1 - \omega)\, d\omega X_k^2\right)\right.\right.
\\&\left.\left.
\times \left(n^{-1/2}(Z_k-  Y)' f_{2n}'(W_{kn}) + \frac{1}{2}n^{-1}(Z'_k f''_{2n}(W_{kn}) Z_k
- Y' f''_{2n}(W_{kn}) Y)\right. \right.\right.\\&\left.\left.\left. +
\frac{1}{6}n^{-3/2} \int_0^1 (1-\omega)^2(f_{2n}'''(W_{kn}+ \omega Z_k) Z_k^3-
f_{2n}'''(W_{kn}+ \omega Y) Y^3)\,d\omega\right)\right]\right|.
\end{split}
\end{equation}
After opening the brackets in the right-hand side of (\ref {4.21}) it remains to estimate each of
  the resulting addendums independently. The estimates are performed in the same way,
 using (\ref{3.55}, \ref{3.56}, \ref{3.57}, \ref{4.12} - \ref{4.16}). Therefore,
 we estimate only three of them.

Using (\ref{4.15}), we get
\begin{equation}\label{4.22}
\begin{split}&
\left|n^{-3/2}E\left[\int_0^1 f''_{1n}(S_{knX} + \omega X_k) (1 -
\omega_1)\, d\omega_1 X_k^2\right.\right.\\& \left.\left.\times
\int_0^1 (1-\omega)^2(f_{2n}'''(W_{kn}+ \omega_2 Z_k) Z_k^3-
f_{2n}'''(W_{kn}+ \omega_2 Y) Y^3)\, d\omega_2\right]\right|\\& \le
Cn^{-3/2}c_{n}^{-3}\epsilon_{1n}^{-2} b_n^{\lambda-1}
G_{kn}(\gamma_n) \le C\epsilon_{1n}^{-2}c_{n1}^{-3}
b_n^{2+\lambda}G_{kn}(\gamma_n).
\end{split}
\end{equation}
The first inequality in (\ref{4.22}) is obtained on the base of the following reasoning
\begin{equation}\label{4.23}
\begin{split}&
W_{kn} + n^{-1/2} Z  \in n^{1/2}b_nv + U_{c_n}  \Rightarrow W_{kn}
\in n^{1/2}b_nv + U_{\epsilon n^{-1/2}b_n^{-1}+c_n}\\&
\Rightarrow  n^{1/2}(n-1)^{-1/2} W_{kn}\in (n-1)^{1/2}b_{n-1}v +
(n(n-1)^{-1/2}b_{n} - (n-1)^{1/2}b_{n-1})v \\& + n^{1/2}
(n-1)^{-1/2} U_{\epsilon n^{-1/2}b_n^{-1} +c_n}\\& \Rightarrow
n^{1/2}(n-1)^{-1/2} W_{kn}\in (n-1)^{1/2}b_{n-1}v +
U_{\gamma_n}.
\end{split}
\end{equation}
Using (\ref{4.12}), we get
\begin{equation}\label{4.24}
\begin{split}&
 E[|f'_{1n}(S_{k,n-1,X}) X_k n^{-1} f_{2n}''(W_{kn})
Z_k^2|]\\&\le Cn^{-1}b_n^\lambda c_n^{-2}\epsilon_{1n}^{-1}
G_{kn}(\gamma_n) \le C b_n^{2+\lambda} \epsilon_{1n}^{-1} c_{n1}^{-2} G_{kn}(\gamma_n).
\end{split}
\end{equation}
Using (\ref{3.57}), we get
\begin{equation}\label{4.25}
\begin{split}&
n^{-1/2} E [f'_{1n}(S_{knX}) X_k (Z_k - Y) f'_{2n}(W_{kn})]\\&
=n^{-1/2} E[X_k Z_k] E[ f'_{1n}(S_{knX}) f'_{2n}(W_{kn})]  \le
Cn^{-1/2} b_n^{1 + \lambda} \epsilon_{1n}^{-1} c_{n1}^{-1} G_{kn}(\gamma_n).
\end{split}
\end{equation}
This completes the proof of Lemma  \ref{l4.2}.

 We begin the proof of Theorem \ref{t3} with auxilliary estimates.
\begin{equation}\label{4.26}
\begin{split}&
P(Y \in n^{1/2}b_n +   U_{c_{n}}) \le
\exp\{Cc_{n}n^{1/2}b_n\} P(Y \in n^{1/2}b_n + U) \\&
\le a_0 P(Y \in n^{1/2}b_n + U).
\end{split}
\end{equation}
Note that
\begin{equation}\label{4.27}
Y \in (n-1)^{1/2}b_{n-1} v +  U_{\gamma_n} \Rightarrow Y \in
n^{1/2}b_n v +  U_{\omega_n}
\end{equation}
with $\omega_n = \gamma_n + n^{1/2} b_n - (n-1)^{1/2}b_{n-1}$.

Therefore
\begin{equation}\label{4.28}
\begin{split}&
P(Y \in (n-1)^{1/2}b_{n-1} v +  U_{\gamma_n}) \le
P(Y \in n^{1/2}b_n v +  U_{\omega_n})\\& \le
C\exp\{n^{1/2} b_n\omega_n\} P(Y \in n^{1/2}b_n v +  U)
\le  a_1 P(Y \in n^{1/2}b_n v +  U).
\end{split}
\end{equation}
The further reasonings are based on an induction on $n $. We take a sufficiently large $ n = n_0 $ such that
 $ C n_0 \epsilon_{1n_0}^{-2}c_{n_0,1}^{-3} b_{n_0}^{2+\lambda}   < a$ with $aa_0a_1 <1$.
We take $C_{n_0}$ such that
\begin{equation}\label{4.29}
C_{n_0} P(Y \in n_0^{1/2} b_{n_0} +  U) E [f_{1n}(S_{n_0X})]\ge 1.
\end{equation}
Then
\begin{equation}\label{4.30}
E [f_{1n}(S_{n_0X}),S_{n_0Z} \in n_0^{1/2}b_{n_0} v +  U] \le C_{n_0}
P(Y \in n_0^{1/2} b_{n_0} +  U) E [f_{1n}(S_{n_0X})].
\end{equation}
Suppose Theorem \ref{t3} was proved for $n-1\ge n_0$. Let us prove it for $n$. We show
\begin{equation}\label{4.31}
E [f_{1n}(S_{nX}),S_{nZ} \in n^{1/2}b_{n} v +  U] \le C_{n} P(Y \in
n^{1/2} b_{n} +  U) E [f_1(S_{nX})]
\end{equation}
where $C_n = a_0 + C_{n-1}a a_1$.
Then, since $C_n$  form  geometric progression with exponent
 $ a a_0 a_1 < 1$,
Theorem \ref{t3} follows from (\ref{4.31}).

Applying (\ref{4.17}) and the inductive assumption
, we get
\begin{equation}\label{4.32}
\begin{split}&
E [f_{1n}(S_{nX}),S_{nZ} \in n^{1/2}b_{n} v +  U]\le P(Y \in n^{1/2}
b_{n} +  U_{c_{1n}}) E [f_{1n}(S_{nX})] \\& +Cnb_n^{2+\lambda}
c_{n1}^{-3} \epsilon_{1n}^{-2} C_{n-1} E[f_{1n}(S_{nX})] P(Y \in
(n-1)^{1/2} b_{n-1} +  U_{\gamma_{n}})\\& \le (a_0 +C_{n-1}a
a_1) E[f_{1n}(S_{nX})] P(Y \in n^{1/2} b_{n} +  U).
\end{split}
\end{equation}
This implies Theorem \ref{t3}.

{\it Proof of Theorem \ref{t2}}. In the proofs of Theorem \ref{t2} and  Osypov
Theorem \cite{os} the basic reasonings coinside. The difference is only in the preliminary estimates.
On these estimates the basic reasoning are based on.

Denote $\phi(h) = E [\exp\{h'X\}]$. Define random vector $X_h$
having the conjugate distribution
$$
F_h(dx) = F(dx,h) = \phi^{-1}(h) \exp\{h'x\} F(dx).
$$
Denote
$$
m(h) = E_h [X_h] , \quad \sigma(h) = \mbox{Var} [X_h].
$$
For any $v \in R^d$ denote $h(v)$ the solution of the equation
\begin{equation}\label{4.1}
m(h) = v.
\end{equation}

\begin{lemma} \label{l4.1} For all $v, |v| <\epsilon b_n, \epsilon > 0$
there exists the solution  $h(v)$ of equation (\ref{4.1}) and
\begin{equation}\label{4.2}
\phi(h) = 1 + |h|^2/2 + O(|h|^3 b_n^{\lambda-1}),
\end{equation}
\begin{equation}\label{4.3}
m(h) = h + O(|h|^2 b_n^{\lambda -1}),
\end{equation}
\begin{equation}\label{4.4}
h(v) = v + O(|v|^2 b_n^{\lambda -1}),
\end{equation}
\begin{equation}\label{4.5}
\sigma(h) = I(1 + O(|h|^2b_n^{\lambda-1})).
\end{equation}
\end{lemma}
{\it Proof of Lemma \ref{l4.1}}. Expanding in the Taylor series we  get
\begin{equation}\label{4.6}
\phi(h) = 1 + \frac{1}{2}\int (h'x)^2 \, dF(x) + O\left(|h|^3 \int |x|^3 \, dF(x)\right) =
1 + \frac{1}{2}|h|^2 + O(|h|^3  b_n^{\lambda-1}),
\end{equation}
\begin{equation}\label{4.7}
\begin{split}&
m(h) = \phi^{-1}(h) \int x \exp\{h'x\} \, dF(x) \\& =
\int x (h'x) dF(x)( 1 - |h|^2/2 + O(|h|^3  b_n^{\lambda-1})+ O\left(\int x (h'x)^2 dF(x)\right)\\&
 = h + O(|h|^2 + |h|^2  b_n^{\lambda-1}).
\end{split}
\end{equation}
Substituting (\ref{4.7}) in (\ref{4.1}), we get (\ref{4.4}).
Estimating similarly to (\ref{4.7}), we get (\ref{4.5}).

Denote
\begin{equation}\label{4.8}
\Lambda(h,v) = -(h,v) + \ln \phi(h).
\end{equation}
By (\ref{4.2},\ref{4.4}), we get
\begin{equation}\label{4.9}
\ln \phi(h(v)) =  \frac{1}{2}h^2(v)(1 + O(b_n^\lambda)).
\end{equation}
By (\ref{4.5}), we get
\begin{equation}\label{4.10}
\mbox{det}^{-1/2}\sigma(h(v)) = 1 + O(b_n^\lambda).
\end{equation}
By (\ref{4.4}) and (\ref{4.9}) we get
\begin{equation}\label{4.11}
\begin{split}&
\Lambda(h(v),v) = |v|^2(1 + O(|v|b_n^{\lambda-1})) - \frac{1}{2}|v|^2(1 + O(b_n^\lambda))\\&
=\frac{1}{2} |v|^2 + O(|v|^2b_n^\lambda).
\end{split}
\end{equation}
The estimates (\ref{4.2}-\ref{4.5}) and (\ref{4.9}-\ref{4.11}) are the versions of similar estimates in  \cite{os}.
Using these estimates we get Theorem \ref{t2} on the base of the same reasoning as in \cite{os}. This reasoning is omitted

\section{ Proofs of Lemmas \ref{l3.1},\ref{l3.2},\ref{l3.5},\ref{l3.6} and \ref{l3.10}-\ref{l3.12}}

The Lemmas will be proved in the following order: \ref{l3.1},\ref{l3.2},\ref{l3.5},\ref{l3.6},\ref{l3.8},\ref{l3.10},\ref{l3.9},\ref{l3.11},\ref{l3.12}.

{\it Proof of Lemma \ref{l3.1}.} Let $h \in \Psi_j(\theta)$ and $h_1 \in \Pi(h)$. By (\ref{2.1}) and (\ref{2.3}), we get
\begin{equation}\label{3.13}
\begin{split}&
P_{h_1}(|\eta(h_1,h)| > \epsilon) \le P_{h_1}(|\eta(h_1,h) - \frac {1}{2}
\bar h'\tau_{h_1}| > \epsilon/2) + P_{h_1}(|\bar h'\tau_{h_1}| >
\epsilon/2) \\& < 4\epsilon^{-2} E_{h_1}[(\eta(h_1,h) - \frac {1}{2}
\bar h'\tau_{h_1})^2] +
2^{2+\lambda}\epsilon^{-2-\lambda}|\bar h|^{2+\lambda}
E_{h_1}|\tau_{h_1}|^{2+\lambda} \le
C|\bar h|^{2+\lambda}.
\end{split}
\end{equation}
By straightforward calculations, using (\ref{3.13}), for $1 \le j \le m$, we get
\begin{equation}\label{3.14}
P(V_{h}(\theta)) \le CP_{h_1}(|\eta(h_{1},h)| > \epsilon j^{-2})  \le C \epsilon^{-2}j^4|\bar h|^{2+\lambda}
\le C
j^4\left(\frac{b_{n}}{2^j}\right)^{2
+\lambda}.
\end{equation}
In the case of $j = m+1$ the constant $C$ in (\ref{3.14}) is replaced with $C c_{3n}^{d-1}$.
By (\ref{3.14}), we get
\begin{equation}\label{3.15}
P(B_{4n}(\theta)) < Cn\sum_{j=1}^m
2^{j}\left(\frac{b_n}{2^j}\right)^{2 +\lambda} j^4 + C nc_{3n}^{d-1} 2^m
c_{3n}^{2+\lambda}\delta_{1n}^{2+\lambda}m^4.
\end{equation}
Note that
$2^m =Cc_{1n}^{-1}nb_n^2(1 + o(1)) $ .
Therefore, using $n^{-\lambda}b_n^{-\lambda} < nb_n^{2+\lambda}$, we get
\begin{equation}\label{3.16}
\begin{split}&
P(B_{4n}(\theta)) < Cnb_n^{2+\lambda} + Cn
c_{3n}^{d+1+\lambda} 2^{-m(1+\lambda)} m^4
b_n^{2+\lambda}
\\& \le Cnb_n^{2+\lambda}\epsilon^{-2-\lambda} + C  C_n c_{3n}^{d + 2 +\lambda}
n^{-\lambda}b_n^{-\lambda} m^4 =
O(nb_n^{2+\lambda})= o(1)
\end{split}
\end{equation}
 if $ c_{3n}$ tends to infinity sufficiently slowly.

Since $P^{(s)}_{h,h_1}(S) < C |\bar h|^{2+\lambda} $,
then, arguing similarly (\ref{3.14})-(\ref{3.16}), we get
\begin{equation}\label{3.16a}
\begin{split}&
P(D_{nile}) \le Cn\sum_{j=1}^{m+1}\sum_{h\in \Psi_j(\theta)} P^{(s)}_{h,h_1}(S)\\&
\le Cn\sum_{j=1}^m 2^j (b_n2^{-j})^{2+\lambda} + Cn c_{3n}^{d+1+\lambda}2^m
\delta_{1n}^{2+\lambda} = o(1).
\end{split}
\end{equation}
Now (\ref{3.16},\ref{3.16a}) implies (\ref{3.10}).

{\it Proof of Lemma \ref{l3.2}}. Applying the Chebyshev inequality
and using
 (\ref{2.3}), we get
\begin{equation}\label{3.18}
 P( B_{3n1}) \le \epsilon^{-2-\lambda} b_n^{2+\lambda} E[|\tau|^{2+\lambda}] < C b_n^{2+\lambda}.
\end{equation}
Let $h \in \Psi_j(\theta), 1 \le j \le m+1$. By Chebyshev inequality, we get
\begin{equation}\label{3.18a}
\begin{split}&
P(|\tau_{sh} -\tau_s| > \epsilon b_n^{-1}2^{j/2}|A_{4n1}) <
C 2^{-j(2+\lambda)/2}b_n^{2+\lambda}\epsilon^{-2-\lambda} (E[ |\tau_h|^{2+\lambda}|A_{4n1}] + E[|\tau|^{2+\lambda}])\\&
<C 2^{-j(2+\lambda)/2}b_n^{2+\lambda}\epsilon^{-2-\lambda} (E_h[ |\tau_h|^{2+\lambda}]+ E[|\tau|^{2+\lambda}]) \le
C 2^{-j(2+\lambda)/2}b_n^{2+\lambda}.
\end{split}
\end{equation}
By (\ref{3.18}), (\ref{3.18a}), we get
\begin{equation}\label{3.18b}
P(B_{3nile}) < Cn\sum_{j=1}^m 2^j b_n^{2+\lambda}  2^{-j(2+\lambda)/2}+ Cnc_{3n}^{d-1}2^m
2^{-m(2+\lambda)/2}b_n^{2+\lambda} < C  n  b_n^{2+\lambda} = o(1).
\end{equation}
By (\ref{3.16}),(\ref{3.16a}) and (\ref{3.18b}), we get
\begin{equation}\label{3.18c}
P(B_{1nile}) < C  n b_n^{2+\lambda}.
\end{equation}

{\it Proof of Lemma \ref{l3.5}} Since $E[\tau] =0$, we have
\begin{equation}\label{3.43}
\begin{split}&
|E[\tau, A_{1n1}]|=|E[\tau, B_{1n1}]|\\& \le
E[|\tau|, |\tau| >  b_n^{-1}] +
E[|\tau|, B_{1n1}\cap \{|\tau| \le  b_n^{-1}\})]\\&
\le b_n^{1+\lambda} E|\tau|^{2+\lambda} +
 b_n^{-1}  P(B_{1n1}) =
O(b_n^{1+\lambda})
\end{split}
\end{equation}
where the last equality follows from (\ref{2.3}),(\ref{3.16}),(\ref{3.18}).

The proof of (\ref{3.42}) is similar and is omitted.

The considerable part of the subsequent estimates is based on the following lemma.
\begin{lemma} \label{l61} Let  $h \in \Psi_j(\theta), h_1 \in \Pi(h), 1 \le j \le m+1, \theta \in \Theta_{nile}$. Then, for any $a \ge 0, b \ge 0, a + b \ge 2 + \lambda$, there
holds
\begin{equation}\label{71}
E_{h_1}[|\bar h\tau_{h_1}|^a |\eta(h_1,h)|^b, A_{1n1}] \le C |\bar h|^{2+\lambda}.
\end{equation}
\end{lemma}
{\it Proof of Lemma \ref{l61}}. By (\ref{2.1}) and (\ref{2.3}, we get
\begin{equation}\label{72}
\begin{split}&
E_{h_1}[|\bar h\tau_{h_1}|^a |\eta(h_1,h)|^b, A_{1n1}] \le
CE_{h_1}[|\bar h\tau_{h_1}|^{a+b}, A_{1n1}] +
CE_{h_1}[|\eta(h_1,h)|^{a+b}, A_{1n1}]\\& \le
CE_{h_1}[|\bar h\tau_{h_1}|^{a+b}, A_{1n1}] +
CE_{h_1}[|\eta(h_1,h) - \bar h\tau_{h_1}|^{a+b}, A_{1n1}]\\& \le
CE_{h_1}[|\bar h\tau_{h_1}|^{2+\lambda}, A_{1n1}] +
CE_{h_1}[|\eta(h_1,h) - \bar h\tau_{h_1}|^{2}, A_{1n1}]
\le C |\bar h|^{2+\lambda}.
\end{split}
\end{equation}

{\it Proof of Lemma \ref{l3.6}}. Expanding $\xi_{n}$
in the Taylor series, we get
\begin{equation} \label{5.1}
S_{n\theta} = \sum_{s=1}^n (2\eta_{ns}(\theta) - \theta' \tau_{s})
- \sum_{s=1}^n \eta^2_{ns}(\theta)  + \frac{2}{3}\sum_{s=1}^n
\frac{\eta_{ns}^3(\theta)}{(1 + \kappa\eta_{ns}(\theta))^3} +
2n\rho^2(0,\theta)
\end{equation}
where $0 \le \kappa \le 1$.

Since $ E [\eta_n^2(\theta)] = \rho^2(0,\theta) $ and $ 2 E
 [\eta_n (\theta)] =-E [\eta_n^2 (\theta)] = - \rho^2(0,\theta) $,
 by virtue of  (\ref{2.2}), we get
\begin{equation} \label{5.2}
E[(2\eta_{n}(\theta) - \theta' \tau) -
\eta^2_{ns}(\theta) + \frac{1}{2} \theta' I\theta] = O(|\theta|^{2+\lambda}).
\end{equation}
By (\ref{3.16},\ref{3.18c}), we get
\begin{equation} \label{5.3}
\begin{split}&
E[|\eta_n(\theta)| , B_{1n1}) \le E[|\eta_n(\theta)| , |\eta_n(\theta)| > \epsilon] +
E[|\eta_n(\theta)| , B_{1n1} \setminus \{|\eta_n(\theta)| < \epsilon\}]\\&
\le E[|\eta_n(\theta)| , |\eta_n(\theta)| > \epsilon ] +   \epsilon P(B_{1n1}) \le E[|\eta_n(\theta)|,
|\eta_n(\theta)| > \epsilon ] +
C b_n^{2+\lambda}.
\end{split}
\end{equation}
By (\ref{2.1}, \ref{2.3}), we get
\begin{equation} \label{5.4}
\begin{split}&
E[|\eta_n(\theta)| , |\eta_n(\theta)| > \epsilon ] \\& \le
E[|\eta_n(\theta)| , |\eta_n(\theta)| > \epsilon, |\eta_n(\theta) - \frac {1}{2}\theta' \tau| < \epsilon/2 ]  +
E[|\eta_n(\theta)| , |\eta_n(\theta)| > \epsilon, |\theta \tau| < \epsilon/2 ]\\&
\le C E [|\theta'\tau|, |\eta_n(\theta)| > \epsilon, |\eta_n(\theta) - \frac {1}{2}\theta' \tau| < \epsilon/2]
+ 4 \epsilon^{-1}E[(\eta_n(\theta) - \frac {1}{2}\theta' \tau)^2] \\& \le
C\epsilon^{-1-\lambda} E[|\theta'\tau|^{2 + \lambda}] +
Cb_n^{2+\lambda} \le C b_n^{2+\lambda}.
\end{split}
\end{equation}
By (\ref{5.3}) and (\ref{5.4}), we get
\begin{equation}\label{5.4a}
E[\eta_n(\theta)|B_{1n1}] \le C  b_n^{2+\lambda}.
\end{equation}

Arguing similarly to (\ref{5.3},
\ref{5.4}), we get
\begin{equation}\label{5.5}
E[\eta_n^2(\theta) , B_{1n1}] = O( b_n^{2+\lambda}).
\end{equation}
By (\ref{5.2},\ref{3.18c},\ref{3.43},\ref{5.4a}),(\ref{5.5}), we get
\begin{equation} \label{5.6}
E[(2\eta_{n}(\theta) - \frac {1}{2}\theta'\tau) -  \eta^2_{ns}(\theta) +
\frac{1}{2} \theta' I\theta, B_{1n1}] =
O(|b_n|^{2+\lambda}).
\end{equation}
By Lemma \ref{l61}, we get
\begin{equation} \label{5.7}
E\left[\left|\frac{\eta_n^3(\theta)}{(1 +
\kappa\eta_n(\theta))^3}\right|, A_{1n1}\right]  \le
C E[|\eta_n^3(\theta)|, A_{1n1}]  \le
  C|\theta|^{2+\lambda}.
\end{equation}
By (\ref{5.1}),(\ref{5.2}),(\ref{5.6}),(\ref{5.7}) we get (\ref{3.20}).

{\it Proof of Lemma \ref{l3.8}}. Using (\ref{5.1}), we get
 \begin{equation}\label{5.8}
\begin{split} &
E[(\xi(\theta) - \theta'\tau)^2, A_{1n1}] \le C
E[(\eta_n(\theta) - \frac {1}{2}\theta' \tau)^2] \\&
 + C E[\eta_n^4(\theta), A_{1n1}]  + C E[\eta_n^6(\theta), A_{1n1}].
\end{split}
\end{equation}
By Lemma \ref{l61}, we get
\begin{equation}\label{5.9a}
E[\eta_n^4(\theta), A_{1n1}]
 = O(|\theta|^{2+\lambda}).
\end{equation}
and
\begin{equation}\label{5.10}
E[\eta_n^6(\theta), A_{1n1}]= O(|\theta|^{2+\lambda}).
\end{equation}
By (\ref{2.1}), (\ref{5.8}), (\ref{5.9a}), (\ref{5.10}) we get (\ref{3.61}).

Estimating similarly to (\ref{5.8}-\ref{5.10}), we get
\begin{equation}\label{5.11}
\begin{split} &
E[(\xi(h_1,h) - \frac {1}{2}\bar h' \tau_{h_1})^2, A_{1n1}] \\& \le C
E_{h_1}[(\xi(h_1,h) - \frac {1}{2}\bar h' \tau_{h_1})^2, A_{1n1}] \le C |\bar
h|^{2 +\lambda}.
\end{split}
\end{equation}
This implies (\ref{3.62}).

{ \it Proof of Lemma \ref{l3.10}}. Applying the Cauchy inequality, by (\ref{3.63}), we get
\begin{equation}\label{5.15}
\begin{split} &
E[(\xi(h_1,h) - \bar h'\tau_{h_1})(v'\tau), A_{1n1}] \\&
\le (E[(\xi(h_1,h) - \bar h'\tau_{h_1})^2, A_{1n1}])^{1/2}(E[(v'\tau)^2, A_{1n1}])^{1/2}
\le C|v|  |\bar h|^{1+\lambda/2}.
\end{split}
\end{equation}
This completes the proof of Lemma \ref{l3.10}.

{\it Proof of Lemma \ref{l3.9}.}
 Using the inequality
$
(a+b)^2 -2b^2 \le 2a^2,
$
putting $a =\eta(0,u)+ \frac {1}{2} u'\tau- \eta(h,h+u)+\frac {1}{2}u'\tau_h $ and $b = \eta(h,h+u)-\eta(0,u)$, we get
\begin{equation}\label{5.12a}
\begin{split}&
E[(u'(\tau - \tau_h))^2, A_{1n1}]- 2E[(\eta(h,h+u)
-\eta(0,u))^2, A_{1n1}] \\& \le 2
E[(\eta(h,h+u)-\frac {1}{2}u'\tau_h -\eta(0,u)+ \frac {1}{2}u'\tau)^2,  A_{1n1}]\doteq J.
\end{split}
\end{equation}
Using the inequality  $2a^2 \le 4(a+b)^2 +4b^2$,  putting $a = \eta(h,h+u)-\frac {1}{2}u'\tau_h -\eta(0,u)+ \frac {1}{2}u'\tau$
and $b = \eta(0,u) - \frac {1}{2}u'\tau$, by (\ref{2.1}), we get
\begin{equation}\label{5.13}
\begin{split}&
J  \le 4 E[(\eta(h,h+u) -\frac {1}{2}u'\tau_h)^2,  A_{1n1}] +4E[(\eta(0,u)- \frac {1}{2}u'\tau)^2,  A_{1n1}]\\&   \le
CE_h[(\eta(h,h+u) -\frac {1}{2}u'\tau_h)^2] + C|u|^{2+\lambda} \le C|u|^{2+\lambda}.
\end{split}
\end{equation}
Thus, for the proof of  (\ref{3.63}), it suffices to show
\begin{equation}\label{5.14}
J_1 \doteq E[(\eta(h,h+u) -\eta(0,u))^2,  A_{1n1}] = O(|u|^2 |h|^\lambda).
\end{equation}
By straightforward calculations, we get
$$
(\eta(h,h+u) - \eta(0,u))^2
$$
$$
=(\eta(0,h+u) -\eta(0,h) - \eta(0,u) -
\eta(0,h)\eta(0,u))^2(\eta(0,h)+1)^{-2}.
$$
Therefore
\begin{equation}\label{5.14a}
\begin{split}&
J_1 = E[(\eta(0,h+u) -\eta(0,h) - \eta(0,u) -
\eta(0,h)\eta(0,u))^2(\eta(0,h)+1)^{-2},  A_{1n1}]\\&\le C E[(\eta(0,h+u) -\eta(0,h) - \eta(0,u)
- \eta(0,h)\eta(0,u))^2,  A_{1n1}]\\&\le
C E[(\eta(0,h+u)- \frac {1}{2}(h+u)'\tau -(\eta(0,h) - \frac {1}{2}h'\tau) -
(\eta(0,u) - \frac {1}{2}u'\tau))^2,  A_{1n1}] \\& +
C E[\eta^2(0,h)\eta^2(0,u)), A_{1n1}] \doteq J_{11} +
J_{12}.
\end{split}
\end{equation}
Applying (\ref{2.1}), we get
\begin{equation}\label{5.14b}
\begin{split}&
J_{11} \le C E[(\eta(0,h+u)- \frac{1}{2}(h+u)'\tau)^2]  + C
E[ (\eta(0,h) - \frac{1}{2}h'\tau)^2 ] \\&   + C
E[(\eta(0,u) - \frac {1}{2}u'\tau)^2] \le C |h+ u|^{2+\lambda} + C |h|^{2+\lambda}.
\end{split}
 \end{equation}
By Lemma \ref{l61}, we get
\begin{equation}\label{5.14c}
J_{12} \le C E[\eta^4(0,h),
 A_{1n1}] +
 C E[\eta^4(0,u),  A_{1n1}]
\le C (|u|^{2+\lambda} + |h|^{2+\lambda} ).
 \end{equation}
By (\ref{5.14a}-\ref{5.14c},\ref{5.13},\ref{5.12a}), we get
\begin{equation}\label{5.14d}
E[(u'(\tau-\frac {1}{2}\tau_h))^2,  A_{1n1}] \le
C(|h+u|^{2+\lambda} + |u|^{2+\lambda} + |h|^{2+\lambda}
).
\end{equation}
Putting $ |u| = c_0 |h|$ and  $C_1 = C((1+ c_0)^{2+\lambda} + c_0^{2+\lambda} + c_0^2)c_0^{-2}$, we get
\begin{equation}\label{5.14e}
E[(u'(\tau - \tau_h))^2, A_{1n1}] \le C_1|u|^2
|h|^\lambda.
\end{equation}
This completes the proof of Lemma \ref{l3.9}.

{\it Proof of Lemma \ref{l3.11}}. Denote
\begin{equation}\label{5.22}
\begin{split}&
W \doteq E[(h'_1\tau)(\xi(h_{1},h) - \bar h' \tau_{h_{1}}
)|A_{1n1}] = E[(h'_1(\tau - \tau_{h_1}))(\xi(h_{1},h) - \bar h' \tau_{h_{1}}
)|A_{1n1}]\\& +
E[(h'_1\tau_{h_1})(\xi(h_{1},h) - \bar h' \tau_{h_{1}}
)|A_{1n1}] \doteq W_{11} + W_{12}.
\end{split}
\end{equation}
By (\ref{3.63}),(\ref{3.62}), we get
\begin{equation}\label{80}
\begin{split}&
W_{11} \le (E[(h'_1(\tau - \tau_{h_1}))^2|A_{1n1}])^{1/2}
(E[(\xi(h_{1},h) - \bar h' \tau_{h_{1}}
)^2|A_{1n1}])^{1/2}\\& \le C|h_1|^{1+\lambda/2} |\bar h|^{1 + \lambda/2}.
\end{split}
\end{equation}
We have
\begin{equation}\label{81}
\begin{split}&
W_{12} = E_{h_1}[(1 + \eta(h_1,0))^2(h'_1\tau_{h_1})(\xi(h_{1},h) - \bar h' \tau_{h_{1}}
)|A_{1n1}]\\&
= E_{h_1}[(h'_1\tau_{h_1})(\xi(h_{1},h) - \bar h' \tau_{h_{1}}
)|A_{1n1}] +2 E_{h_1}[ \eta(h_1,0)(h'_1\tau_{h_1})(\xi(h_{1},h) - \bar h' \tau_{h_{1}}
)|A_{1n1}]\\& + E_{h_1}[ \eta^2(h_1,0)(h'_1\tau_{h_1})(\xi(h_{1},h) - \bar h' \tau_{h_{1}}
)|A_{1n1}] \doteq W_{121} + W_{122} + W_{123}.
\end{split}
\end{equation}
By (\ref{5.1}), we get
\begin{equation}\label{82}
\begin{split}&
W_{121} = E_{h_1}[h'_1\tau_{h_1}(2\eta(h_1,h) - \bar h \tau_{h_1}), A_{1n1}] -
 E_{h_1}[h'_1\tau_{h_1}\eta^2(h_1,h), A_{1n1}] \\& +
\frac{2}{3}E_{h_1}\left[h'_1\tau_{h_1}\frac{\eta^3(h_1,h)}{(1 + \kappa\eta(h_1,h))^3}, A_{1n1}\right] \doteq
W_{1211} + W_{1212} + W_{1213}.
\end{split}
\end{equation}
By (\ref{2.1}),(\ref{2.2}), we get
\begin{equation}\label{5.34}
\begin{split}&
O(|\bar h|^{2+\lambda}) =E_{h_1}[(\eta(h_1,h) - \frac {1}{2}\bar h'
\tau_{h_1})^2] = \rho^2(h_1,h) - E_{h_1}[\eta(h_1,h)\bar h'\tau_{h_1}] +\frac{1}{4}\bar h
I(h_1) \bar h\\&
= \frac{1}{2} \bar h' I(h_1) \bar h (1 +|\bar h|^\lambda) -
E_{h_1}[\eta(h_1,h)\bar h\tau_{h_1}].
\end{split}
\end{equation}
Since $h_1 \parallel \bar h$, by (\ref{5.34}), we get
\begin{equation}\label{5.35}
E_{h_1}[h_1'\tau_{h_1}\eta(h_1,h)] =  \frac{1}{2}h_1' I(h_1) \bar h (1 +O(|\bar
h|^\lambda)).
\end{equation}
Applying the Holder's inequality, we get
\begin{equation}\label{5.36}
\begin{split}&
E_{h_{1}}[h'_1\tau_{h_{1}}(\eta(h_{1},h) - \frac {1}{2}\bar h'
\tau_{h_1}), B_{1n1}]\\& \le (E_{h_1}[(h_1'\tau_{h_1})^{2+\lambda}])^{\frac{1}{2+\lambda}}(E_{h_1}
[(\eta(h_1,h) - \frac {1}{2}\bar h \tau_{h_1})^2])^{1/2}
(P_{h_1}(B_{1n1}))^{\frac{\lambda}{2(2+\lambda)}}\\& = O(|h_1| |\bar h|^{1+\lambda/2}
b_n^{\lambda/2}).
\end{split}
\end{equation}
By (\ref{5.35}),(\ref{5.36}),(\ref{3.42}), we get
\begin{equation}\label{84}
W_{1211} = O(|h'_1| |\bar h| b_n^\lambda).
\end{equation}
By Lemma \ref{l61}, we get
\begin{equation}\label{85}
W_{1212} + W_{1213} = O(|h_1| |\bar h|^{1+\lambda}.
\end{equation}
By (\ref{82}),(\ref{84}),(\ref{85}), we get
\begin{equation}\label{86}
W_{121}= O(|h'_1| |\bar h| b_n^\lambda).
\end{equation}
Using Lemma \ref{l61} and (\ref{5.1}), we get
\begin{equation}\label{87}
W_{122} + W_{123} = O(|\bar h|^{1+\lambda} |h_1|).
\end{equation}
By (\ref{81}), (\ref{86}), (\ref{87}), we get
\begin{equation}\label{88}
W_{12} = O(|h'_1| |\bar h| b_n^\lambda).
\end{equation}
By (\ref{5.22}), (\ref{80}), (\ref{88}), we get (\ref{3.65}).

{\it Proof of Lemma \ref{l3.12}}. We begin with the proof of (\ref{3.63c}).
Using (\ref{3.63}), we get
\begin{equation}\label{6.1}]
E[\bar h'(\tau - \tau_{h_1})\tau_k,  A_{1n1}] \le
(E[\bar h'(\tau - \tau_{h_1})^2,  A_{1n1}])^{1/2} (E[\tau_k^2])^{1/2} < C |\bar h| |h_1|^{\lambda/2}.
\end{equation}
The proof of (\ref{3.63d}) is based on the following reasoning. By (\ref{3.63}), we get
\begin{equation}\label{73}
\begin{split}&
O(|\bar h|^2 b_n^{\lambda}) = E[(\bar h(\tau - \tau_{h_1}))^2,A_{1n1}] = E[(\bar h\tau)^2, A_{1n1}] -\\&
-2E[(\bar h\tau)(\bar h\tau_h), A_{1n1}]  + E[(\bar h\tau_{h_1})^2, A_{1n1}] \doteq J_1 - 2J_2 + J_3.
\end{split}
\end{equation}
We have
\begin{equation}\label{74}
\begin{split}&
J_3 = E_{h_1} [(\eta(h_1,0) +1)^2(\bar h\tau_{h_1})^2, A_{1n1}]\\& =
E_{h_1} [\eta^2(h_1,0) (\bar h\tau_{h_1})^2, A_{1n1}]
+ 2 E_{h_1} [\eta(h_1,0) (\bar h\tau_{h_1})^2, A_{1n1}]\\&
+ E_{h_1}[(\bar h\tau_{h_1})^2, A_{1n1}] = J_{31} + 2J_{32} + J_{33}.
\end{split}
\end{equation}
By Lemma \ref{l61}, we get
\begin{equation}\label{75}
J_{31} + 2J_{32} \le C |\bar h|^2 |h|^\lambda.
\end{equation}
Estimating similarly to the proof of (\ref{3.41}),(\ref{3.42}), we get
\begin{equation}\label{76}
J_{33} = \bar h' I(h) \bar h + O(|\bar h|^2 b_n^\lambda).
\end{equation}
By (\ref{74})-(\ref{76}), we get
\begin{equation}\label{77}
J_3 = \bar h_1' I(h_1) \bar h_1 + O(|\bar h|^2 b_n^\lambda).
\end{equation}
By (\ref{73}), (\ref{3.42}),(\ref{77}),  we get
\begin{equation}\label{78}
J_2 = \bar h_1' I \bar h_1 + O(|\bar h|^2 b_n^\lambda).
\end{equation}
By (\ref{78}),(\ref{3.42}), we get
\begin{equation}\label{79}
J_1 - J_2 = O(|\bar h|^2b_n^\lambda).
\end{equation}
This implies (\ref{3.63d}).

\bigskip
\end{document}